\documentclass[twocolumn]{autart}     

\usepackage{amssymb,amsmath,amsopn}
\usepackage{amsfonts}
\usepackage{amstext}
\usepackage{empheq}
\usepackage{epsfig}
\usepackage{color}
\usepackage{graphicx,colordvi}
\usepackage{epstopdf}
\usepackage{float}
\usepackage{CJK}

\newtheorem{lemma}{ Lemma}[section]
\newtheorem{theorem}{ Theorem}[section]
\newtheorem{definition}{ Definition}[section]

\newtheorem{remark}{ Remark}[section]
\newtheorem{example}{ Example }
\newtheorem{corollary}{ Corollary}
\newtheorem{proposition}{ Proposition}[section]

\allowdisplaybreaks[4]
\begin{document}

\begin{frontmatter}

\title{Mean Field Linear Quadratic Control: Uniform Stabilization and Social Optimality\thanksref{footnoteinfo}}

\thanks[footnoteinfo]{This work was supported by National Key R\&D Program of China under Grant 2018YFA0703800, and the National Natural Science Foundation of China under Grants 61573221, 61633014, 61773241 and 61877057.}

\author[sdu]{Bing-Chang Wang}\ead{bcwang@sdu.edu.cn},
\author[sdu]{Huanshui Zhang}	\ead{hszhang@sdu.edu.cn},
\author[amss,ams]{Ji-Feng Zhang}\ead{jif@iss.ac.cn}

\address[sdu]{School of Control Science and Engineering,
Shandong University, Jinan 250061, P. R. China.}	

\address[amss]{the Key Laboratory of Systems and Control, Academy of Mathematics and Systems Science, Chinese Academy of Sciences, Beijing 100190, China}

\address[ams]{School of Mathematical Sciences, University of Chinese Academy of Sciences, Beijing 100149, China }

\begin{abstract}
This paper is concerned with uniform stabilization and social optimality for general mean field linear quadratic control systems,
where subsystems are coupled via individual dynamics and costs, and the state weight is \emph{not assumed with the definiteness condition}.
For the finite-horizon problem, we first obtain a set of forward-backward stochastic differential equations (FBSDEs) from variational analysis, and construct a feedback-type control by decoupling the FBSDEs. For the infinite-horizon problem, by using solutions to two Riccati equations, we design a set of decentralized control laws, which is further proved to be asymptotically social optimal. Some equivalent conditions are given for uniform stabilization of the systems in different cases, respectively.
Finally, the proposed decentralized controls are compared to the asymptotic optimal strategies in previous works. 

\end{abstract}


\end{frontmatter}

%

\section{Introduction}
Mean field games 
have drawn increasing attention in many fields including system control, applied mathematics and economics \cite{BFY13}, \cite{C14}, \cite{GS13}. The mean field game involves a very large population of small interacting players with the feature that while the influence of each one is negligible, the impact of the overall population is significant. By combining mean field approximations and individual's best response,
the dimensionality difficulty is overcome. Mean field games and control
have found wide applications, including smart grids \cite{MCH13}, \cite{CBM15}, \cite{LML19}, finance, economics \cite{GLL11}, \cite{CS14}, \cite{WH15}, 
and social sciences 
\cite{BTN16}, etc.


By now, mean field games have been intensively studied in the LQ (linear-quadratic) framework \cite{HCM07}, \cite{LZ08}, \cite{WZ13}, \cite{ELN13}, \cite{BSYY16}, \cite{MB17}.  
Huang \emph{et al.} developed the Nash certainty equivalence (NCE) based on the fixed-point method and designed an $\epsilon$-Nash equilibrium for mean field LQ games with discount costs by the NCE approach 
\cite{HCM07}. The NCE approach was then applied to the cases with long run average costs \cite{LZ08} and with Markov jump parameters \cite{WZ13}, respectively. The works \cite{CD13}, \cite{BSYY16} employed the adjoint equation approach and the fixed-point
theorem to obtain sufficient conditions for the existence of the
equilibrium strategy over a finite horizon.
For other aspects of mean field games, readers are referred to  \cite{HMC06}, \cite{LL07}, \cite{YMMS12}, \cite{CD13} for nonlinear mean field games, \cite{weintraub2008markov} for oblivious equilibrium in dynamic games, \cite{H10}, \cite{WZ12} for mean field games with major players, \cite{HH16}, \cite{MB17} for robust mean field games.

Besides noncooperative games, social optima in mean field models have also attracted much interest. The social optimum control refers to that all the players cooperate to optimize the common social cost---the sum of individual costs, which is 
a type of team decision problem 
\cite{H80}. Huang \emph{et al.} considered social optima in mean field LQ control, and provided an asymptotic team-optimal solution \cite{HCM12}. Wang and Zhang \cite{WZ17} investigated the mean field social optimal problem where the Markov jump parameter appears as a common source of randomness. 
For further literature, see \cite{HN16} for social optima in mixed games, \cite{AM15} for team-optimal control with finite population and partial information. 

{Most previous results on mean field games and control were given by using the fixed-point method \cite{HCM07}, \cite{LZ08}, \cite{WZ12}, \cite{HCM12}, \cite{C12}, \cite{CD13}, \cite{WZ17}. However, 
the fixed-point analysis (e.g., from the contraction mapping theorem) is sometimes conservative, particularly for high-dimensional systems. In this paper, we 
solve the problem by decoupling directly high-dimensional
forward-backward stochastic differential equations (FBSDEs). In recent years, some 
 progress has been made for study of the optimal LQ control by tackling the FBSDEs. See \cite{Y13}, \cite{ZX17}, \cite{ZQ16}, \cite{SLY16} for details.}

This paper investigates uniform stabilization and social optimality for linear quadratic mean field control systems, 
where subsystems (agents) are coupled via dynamics and individual costs. The state weight $Q$ {is not limited to positive semi-definite}. This model can be taken as a generation of robust mean field control problems \cite{HH16}, \cite{MB17}, \cite{WH17}.
Since the weight $Q$ in the cost functional is indefinite, the prior boundedness of the state is not implied directly by the finiteness of the cost, which brings about additional difficulty to show the social optimality of decentralized control. 

For the finite-horizon social control problem, we first obtain a set of FBSDEs by examining the variation of the social cost, and give a centralized feedback-type control laws by decoupling the FBSDEs. With mean field approximations, we design a set of decentralized control laws. By exploiting the uniform convexity property of the problem, the decentralized controls are further shown to have asymptotic social optimality.
For the infinite-horizon case, we design a set of decentralized control laws by using solutions of two Riccati equations, which is shown to be asymptotically social optimal. Some equivalent conditions are further given for uniform stabilization of all the subsystems when the state weight $Q$ is positive semi-definite or only symmetric.
Furthermore, the explicit expressions of optimal social costs
 are given in terms of the solutions to two Riccati equations, and the proposed decentralized control laws are compared to the feedback strategies in previous works. 
Finally, some numerical examples are given to illustrate the effectiveness of the proposed control laws.

The main contributions of the paper are summarized as follows.
\begin{itemize}
\item
 We first obtain necessary and sufficient existence conditions of finite-horizon centralized optimal control by variational analysis, and then design a feedback-type decentralized control by tackling FBSDEs with mean field approximations.

\item In the case $Q\geq0$, the necessary and sufficient conditions are given for uniform stabilization of the systems with the help of the system's observability and detectability.

\item In the case that $Q$ is \emph{indefinite}, the necessary and sufficient conditions are given for uniform stabilization of the systems using the Hamiltonian matrices.

\item The asymptotically optimal decentralized controls are obtained under very basic assumptions (without verifying the fixed-point condition). The corresponding social costs
 are explicitly given by virtue of the solutions to two Riccati equations.
\end{itemize}


The organization of the paper is as follows. In Section II, the socially optimal control problem is formulated. In Section III, we construct asymptotically optimal  decentralized control laws by tackling FBSDEs for the finite-horizon case. In Section IV, for the infinite-horizon case, the asymptotically optimal controls are designed and analyzed, and some equivalent conditions are further given for uniform stabilization in different cases.
In Section V, some numerical examples are given to show the effectiveness of the proposed control laws. Section VI concludes the paper.

The following notation will be used throughout this paper. $\|\cdot\|$
denotes the Euclidean vector norm or Frobenius matrix norm. For a vector $z$ and a matrix $Q$, $\|z\|_Q^2= z^TQz$, $tr(Q)$ is the trace of the matrix $Q$, and $Q>0$ ($Q\geq0$) means that $Q$ is positive definite (positive semidefinite). For two vectors $x,y$, $\langle x,y\rangle=x^Ty$.
$C([0,T],\mathbb{R}^n)$ is the space of all $\mathbb{R}^n$-valued continuous functions defined on $[0,T]$,
 and $C_{\rho/2}([0,\infty), \mathbb{R}^n)$ is a subspace of $C([0,\infty),\mathbb{R}^n)$ which is given by $\{f|\int_0^{\infty}e^{-\rho t}\|f(t)\|^2dt<\infty  \}.$
$L^2
_{\mathcal  F}(0, T; \mathbb{R}^k)$ is the space of all $\mathcal{F}$-adapted $\mathbb{R}^k$-valued processes $x(\cdot)$ such that
$\mathbb{E}\int_0^T\|x(t)\|^2dt<\infty$.
For convenience of  presentation, we use $C, C_1,C_2,\cdots$ to
denote generic positive constants, which may vary from place to place.

\section{Problem Description}
\label{sec2.3.1}

Consider a large population systems with $N$ agents. Agent $i$ evolves by the following stochastic differential equation:
\begin{equation}\label{eq1}
\begin{aligned}
dx_i(t) = \ & [Ax_i(t)+Bu_i(t)+Gx^{(N)}(t)+f(t)]dt\\
&+\sigma(t) dW_i(t), \quad  1\leq i\leq N,
\end{aligned}
\end{equation}
where $x_i\in
\mathbb{R}^n$ and $u_i\in\mathbb{R}^r$ are the state and input of the $i$th agent.  $x^{(N)}(t)=\frac{1}{N}\sum_{j=1}^Nx_j(t)$, $f, \sigma\in C_{\rho/2}([0,\infty), \mathbb{R}^n)$.
$\{W_i(t),1\leq i\leq N\}$ are a sequence of independent $1$-dimensional Brownian motions on a complete
filtered probability space $(\Omega,
\mathcal F, \{\mathcal F_t\}_{0\leq t\leq T}, \mathbb{P})$.
The cost function of agent $i$ is given by
\begin{equation}\label{eq2}
\begin{aligned}
J_i(u)= \mathbb{E}\int_0^{\infty}e^{-\rho t}\Big\{&\big\|x_i(t)
-\Gamma x^{(N)}(t)-\eta(t)\big\|^2_{Q}\\
&+\|u_i(t)\|^2_{R}\Big\}dt,
\end{aligned}
\end{equation}
where $\rho>0$ and $Q$, $R$ are symmetric matrices with appropriate dimensions. $Q$ is allowed to be \emph{indefinite}. $R>0$, and $\eta\in C_{\rho/2}([0,\infty), \mathbb{R}^n)$.
Denote $u=\{u_1,
\ldots, u_N\}$. 
The decentralized control set is given by $$
\begin{aligned}
{\mathcal  U}_{d,i} =\Big\{u_i\ \big|\ &
u_i(t)\ \hbox{is adapted to}\  \sigma(x_i(s),0\leq s\leq t),\\
&\mathbb{E}\int_0^{\infty}e^{-\rho t}\|u_i(t)\|^2dt<\infty\Big\}.
\end{aligned}
$$ For comparison, define the centralized control sets as
$$\begin{aligned}
  {\mathcal  U}_{c,i} =\Big\{u_i \big|\ &
u_i(t)\ \hbox{is adapted to}\ \sigma\{\bigcup_{i=1}^N{\mathcal F}_t^{i}\}, \\
&\mathbb{E}\int_0^{\infty}e^{-\rho t}\|u_i(t)\|^2dt<\infty\Big\},
\end{aligned}$$
and $
  {\mathcal  U}_{c}\!\! =\!\big\{\!(u_1,\cdots,u_N) \big|\
\!\!u_i\ \hbox{\!belongs to\!}\ {\mathcal  U}_{c,i}, 1\!\leq \!i\!\leq \!N$\big\},
where ${\mathcal F}_t^{i}=\sigma(x_{i}(0),W_i(s),0\leq s\leq t), i=1,\cdots,N$.

In this paper, we mainly study the following problem.

\textbf{(P)}. Seek a set of decentralized control laws to optimize social cost
for the system (\ref{eq1})-(\ref{eq2}), i.e.,
$\inf_{u_i\in {\mathcal  U}_{d,i}}J_{\rm soc},$
where $J_{\rm soc}=\sum_{i=1}^NJ_i(u).$

\begin{remark}
The related results can be extended to the case of multidimensional Brownian motions trivially. 
Here we consider that $\sigma(t)$ is time-varying and satisfies some growth rate.
For convenience of the statement, we assume $W_i$ is scalar and $\sigma\in C_{\rho/2}([0,\infty), \mathbb{R}^n)$.
  For the finite-horizon problem, our results still hold for the case that the matrices $A, B, G,\cdots$ depend on $t$.
\end{remark}

Assume

\textbf{A1)} The initial states of agents $x_i(0), i=1,...,N$ are mutually independent and have the same mathematical expectation. $x_{i}(0)=x_{i0}$, $\mathbb{E}x_i(0)=\bar{x}_0$, $i=1,\cdots,N$. There exists a constant $C_0$ (independent of $N$) such that $\max_{1\leq i \leq N}\mathbb{E}\|x_i(0)\|^2<C_0$.

\section{The finite-horizon problem}\label{sec3}
For the convenience of design, we first consider the following finite-horizon problem.
$$\textbf{(P1)} \inf_{u\in L^2_{{\mathcal  F}}(0, T; \mathbb{R}^{nr})} J_{\rm soc}^{\rm F}(u),
$$
where $J_{\rm soc}^{\rm F}(u)=\sum_{i=1}^NJ_{i}^{\rm F}(u)$ and ${\mathcal F}_t= \sigma\{\bigcup_{i=1}^N{\mathcal F}_t^{i}\}$. Here

\begin{equation}\label{eq3}
\begin{aligned}
J_{i}^{\rm F}(u)=\mathbb{E}\int_0^{T}e^{-\rho t}\Big\{&\big\|x_i(t)
-\Gamma x^{(N)}(t)-\eta(t)\big\|^2_{Q}\\
&+\|u_i(t)\|^2_{R}\Big\}dt.
\end{aligned}
\end{equation}
We first give equivalent conditions for the convexity of (P1).

\begin{proposition}\label{prop1}
 (i) Problem (P1) is convex in $u$ if and only if
for any $u_i\in L^2_{{\mathcal  F}}(0, T; \mathbb{R}^r)$, $i=1,\cdots,N$,
$$\sum_{i=1}^N\mathbb{E}\int_0^Te^{-\rho t}\Big\{\big\|y_i(t)-\Gamma y^{(N)}(t)\big\|^2_{Q}+\|u_i(t)\|^2_{R}\Big\}dt\geq 0, $$
where $y^{(N)}=\sum_{j=1}^Ny_j/N$ and $y_i$ satisfies
\begin{align}\label{eq4aa}
  &dy_i(t)=[A y_i(t)+Gy^{(N)}(t)+Bu_i(t)]dt,\cr
  &   y_i(0)=0, \  i=1,2,\cdots,N.
  \end{align}

(ii) Problem (P1) is uniformly convex in $u$ if and only if
for any $u_i\in L^2_{{\mathcal  F}}(0, T; \mathbb{R}^r)$, there exists $\gamma>0$ such that
$$\begin{aligned}
  &\sum_{i=1}^N\mathbb{E}\int_0^Te^{-\rho t}\Big\{\big\|y_i(t)-\Gamma y^{(N)}(t)\big\|^2_{Q}+\|u_i(t)\|^2_{R}\Big\}dt\cr
  \geq& \gamma\sum_{i=1}^N\mathbb{E}\int_0^Te^{-\rho t}\|u_i(t)\|^2dt.
  \end{aligned}$$
\end{proposition}

\emph{Proof.} Let $x_i$ and $\acute{x}_i$ be the state processes of agent $i$ 
with the control $v$ and $\acute{v}$,
respectively.
 Take any $\lambda_1\in [0, 1]$ and let $\lambda_2= 1-\lambda_1$.
Then
\begin{align*}
&\lambda_1  J_{\rm soc}^{\rm F} (v) +\lambda_2 J_{\rm soc}^{\rm F}  (\acute{v}) - J_{\rm soc}^{\rm F}(\lambda_1 v+\lambda_2 \acute{v})  \cr
 =& \lambda_1\lambda_2 \!\sum_{i=1}^N\mathbb{E}\!\int_0^T\!\big\{ \|x_i(t)\!-\!\acute{x}_i(t) \!-\!\Gamma(x^{(N)}(t)\!-\!\acute{x}^{(N)}(t))\|_Q^2 \cr
& +{ \|v_i(t)-\acute{v}_i(t)\|^2_R}\big\} dt.
\end{align*}
Denote $u= v-\acute{v}$, and $y_i=x_i-\acute{x}_i$. Thus, $y_i$ satisfies
(\ref{eq4aa}). 
By the definition of (uniform) convexity, the lemma follows.  $\hfill \Box$

By examining the variation of $ {J}_{\rm soc}^{\rm F}$, we obtain the necessary and sufficient conditions for
the existence of centralized optimal control of {(P1)}. To simplify the presentation later, we denote by 
$$\left\{
	\begin{aligned}
\Xi& \stackrel{\Delta}{=}\Gamma^TQ+Q\Gamma-\Gamma^TQ\Gamma,\\
 \bar{\eta}&\stackrel{\Delta}{=}Q\eta-\Gamma^T Q\eta.
 \end{aligned}\right.$$
\begin{theorem}\label{thm1}
	Suppose $R>0$. Then
(P1) has a
set of optimal control laws 
if and only if Problem (P1) is convex in $u$ and
	the following equation system admits a set of solutions $(x_i,p_i, \beta_{i}^{j},i,j=1,\cdots,N)$:
	\begin{equation}\label{eq4a}
	\left\{
	\begin{aligned}
	dx_i(t)= &\big(Ax_i(t)\!-\!B{R^{-1}}B^Tp_i(t)\!+\!Gx^{(N)}(t)\!+\!f(t)\big)dt\\
&+\sigma(t) dW_i(t),\\
	dp_i(t)= &-\big[(A-\rho I)^Tp_i(t)+G^Tp^{(N)}(t)+Qx_i(t)\big]dt\\&+\big[\Xi x^{(N)}(t)+\bar{\eta}(t) \big]dt
	+\sum_{j=1}^N\beta_i^j(t)dW_j(t),\\
	x_i(0)=&{x_{i0}},\quad p_i(T)=0,\quad i=1,\cdots,N,
	\end{aligned}\right.
	\end{equation}
	where $p^{(N)}(t)=\frac{1}{N}\sum_{i=1}^Np_i(t)$,
	and
	furthermore the optimal control is given by $\check{u}_i(t)=-{R^{-1}}B^Tp_i(t)$.
\end{theorem}
{\it Proof.}
Suppose that $\check{u}_i=-R^{-1}B^Tp_i,$ where $ {(p_i, \beta_{i}^{j}},i,j=1,\cdots,N)$ is a set of solutions to the equation system
\begin{equation}\label{eq6}
\begin{aligned}
&dp_i(t)=\alpha_i(t)dt+\beta_i^i(t)dW_i(t)+\sum_{j\not =i}\beta_i^j(t)dW_j(t), \\
&p_i(T)=0, \quad i=1,\cdots, N.
\end{aligned}
\end{equation}
Here $\alpha_i$, $i=1,\cdots,N$ are to be determined. 
 Denote by $\check{x}_i$ the state of agent $i$ under the control $\check{u}_i$. For any $u_i\in L^2_{{\mathcal  F}}(0, T; \mathbb{R}^r) $
 and $\theta\in \mathbb{R}\ (\theta\not=0)$, let $u_i^{\theta}=\check{u}_i+\theta u_i$.
 Denote by $ x_i^{\theta}$ the solution of the following perturbed state equation
$$ \begin{aligned}
dx_i^{\theta}(t)=&\Big[Ax_i^{\theta}(t)+B(\check{u}_i(t)+\theta u_i(t))+f(t)\cr&+\frac{G}{N}\sum_{i=1}^Nx^{\theta}_i(t)\Big]dt+\sigma (t)dW_i(t),\cr
x_i^{\theta}(0)=&x_{i0},\ i=1,2,\cdots,N.
 \end{aligned}$$ Let $y_i=(x_i^{\theta}-\check{x}_i)/\theta$. 
It can be verified that
 $y_i$ satisfies (\ref{eq4aa}).
Then by It\^{o}'s formula, for any $i=1,\cdots,N$,
\begin{align*}
0
=&\mathbb{E}[\langle e^{-\rho T} p_i(T),y_i(T)\rangle-\langle p_i(0),y_i(0)\rangle]\cr
=&\mathbb{E}\!\int_0^T\!e^{-\rho t} \big[\langle \alpha_i(t),y_i(t)\rangle+\langle p_i(t),(A-\rho I)y_i(t)\cr
&\qquad\quad+Gy^{(N)}(t)+Bu_i(t)\rangle\big] dt,
\end{align*}
which implies
\begin{align}\label{eq7}
0= &\sum_{i=1}^N \mathbb{E}\int_0^Te^{-\rho t} \big[\langle \alpha_i(t),y_i(t)\rangle+\big\langle p_i(t),(A-\rho I)y_i(t)\cr
&\qquad+Gy^{(N)}(t)+Bu_i(t)\big\rangle\big] dt\cr
= &\sum_{i=1}^N \mathbb{E}\int_0^Te^{-\rho t} \big[\langle \alpha_i(t)+(A-\rho I)^Tp_i(t),y_i(t)\rangle\cr
&\qquad+\langle G^Tp^{(N)}(t),y_i(t)\rangle+\langle B^Tp_i(t),u_i(t)\rangle\big] dt.
\end{align}
From (\ref{eq3}), we have
\begin{equation}\label{eq5a}
\begin{aligned}
 &{J}_{\rm soc}^{\rm F}(\check{u}+\theta u)-{J}_{\rm soc}^{\rm F}(\check{u})=2\theta I_1+{\theta^2}I_2
\end{aligned}
\end{equation}
where $\check{u}=(\check{u}_1,\cdots,\check{u}_N)$, and
\begin{align*}
I_1\stackrel{\Delta}{=}&\sum_{i=1}^N\mathbb{E}\int_0^Te^{-\rho t}  \big[\big\langle Q\big(\check{x}_i(t)-(\Gamma\check{x}^{(N)}(t)+\eta)\big),\cr & y_i(t)-\Gamma y^{(N)}(t)\big\rangle
 +
\langle R \check{u}_i(t),u_i(t)\rangle \big]dt,\cr
I_2\stackrel{\Delta}{=}&\sum_{i=1}^N\mathbb{E}\int_0^T\!\!e^{-\rho t} \big[\big\|y_i(t)
  \! -\!\Gamma y^{(N)}(t)\big\|^2_{Q}
+\|u_i(t)\|^2_{R}\big]dt.
\end{align*}
Note that (suppressing the time $t$)
\begin{align*}
  &\sum_{i=1}^N\mathbb{E}\int_0^Te^{-\rho t} \big\langle Q\big(\check{x}_i-(\Gamma\check{x}^{(N)}+\eta)\big),\Gamma y^{(N)}\big\rangle dt\cr
  =&
\mathbb{ E}\!\int_0^T\!e^{-\rho t} \Big\langle \Gamma^TQ  \sum_{i=1}^N\big(\check{x}_i\!-\!(\Gamma\check{x}^{(N)}\!+\!\eta)\big),\frac{1}{N}  \sum_{j=1}^Ny_j\Big\rangle  dt\cr
 =& \sum_{j=1}^N \mathbb{E}\int_0^Te^{-\rho t} \Big\langle  \frac{\Gamma^TQ}{N} \sum_{i=1}^N\big(\check{x}_i-(\Gamma\check{x}^{(N)}+\eta)\big), y_j\Big\rangle  dt\cr
 =& \sum_{j=1}^N \mathbb{E}\int_0^Te^{-\rho t} \big\langle  {\Gamma^TQ} \big((I-\Gamma)\check{x}^{(N)}-\eta\big), y_j\big\rangle  dt.
\end{align*}
From (\ref{eq7}), one can obtain that
\begin{align}\label{eq10b}
 I_1
=&\sum_{i=1}^N\mathbb{E}\int_0^Te^{-\rho t}\big\langle R\check{u}_i+B^Tp_i,u_i\big\rangle dt\cr
&+\sum_{i=1}^N\mathbb{E}\int_0^Te^{-\rho t}\Big\langle
Q\big(\check{x}_i-(\Gamma\check{x}^{(N)}+\eta)\big)\cr
&-{\Gamma^TQ} \big((I-\Gamma)\check{x}^{(N)}-\eta\big)+\alpha_i\cr
&+(A-\rho I)^Tp_i+G^Tp^{(N)}, y_i\Big\rangle dt.
\end{align}
From (\ref{eq5a}), $\check{u}$ is a minimizer to Problem (P1) if and only if
$I_2\geq0$ and $I_1=0 $.
By Proposition \ref{prop1}, $I_2\geq0$ if and only if (P1) is convex. $I_1=0 $ is equivalent to
\begin{align*}
\alpha_i=&-\big[(A-\rho I)^Tp_i-\Gamma^TQ \big((I-\Gamma)\check{x}^{(N)}-\eta\big)\\
&+Q\big(\check{x}_i-(\Gamma\check{x}^{(N)}+\eta))+G^Tp^{(N)}\big],\cr
\check{u}_i=&-{R^{-1}}B^Tp_i.
\end{align*}
Thus, we have the optimality system (\ref{eq4a}).
This implies that (\ref{eq4a}) admits a solution $(\check{x}_i,\check{p}_i,\check{\beta}_{i}^{j}, i,j=1,\cdots,N)$.

On other hand, if the equation system (\ref{eq4a}) admits a solution $(\check{x}_i,\check{p}_i, \check{\beta}_{i}^{j}, i,j =1,\cdots,N)$.
Let $\check{u}_i=-R^{-1}B^T\check{p}_i$. If (P1) is convex, then $\check{u}$ is a minimizer to Problem (P1).
$\hfill \Box$


It follows from (\ref{eq4a}) that
\begin{equation}\label{eq10}
\left\{
\begin{aligned}
d{x}^{(N)}(t)\!=\!&\big[(A\!\!+\!\!G){x}^{(N)}(t)\!\!-\!\!B{R^{-1}}B^T{p}^{(N)}(t)\!\!+\!\!f(t)\big]dt\\
&+\frac{1}{N}\sum_{i=1}^N\sigma(t)dW_i(t),\\
d{p}^{(N)}(t)=&-\Big[(A+G-\rho I)^T{p}^{(N)}(t)\\
&+(I-\Gamma)^TQ(I-\Gamma){x}^{(N)}(t)-\bar{\eta}(t)\Big]dt\\
&+\frac{1}{N}\sum_{i=1}^N\sum_{j=1}^N{\beta}_i^j(t)dW_j(t),\\
{x}^{(N)}(0)=&\frac{1}{N}\sum_{i=1}^Nx_{i0},\quad {p}^{(N)}(T)=0.
\end{aligned}\right.
\end{equation}
Let $p_i(t)=P(t)x_i(t)+K(t)x^{(N)}(t)+s(t)$, $t\geq 0$. Then by (\ref{eq4a}), (\ref{eq10}) and It\^{o}'s formula (suppressing the time $t$),
$$\begin{aligned}
dp_i=\ &\dot{P}{x}_idt+P\Big[\big(A{x}_i-B{R^{-1}}B^T(Px_i+Kx^{(N)}+s)\\
&+Gx^{(N)}+f\big)dt+\sigma dW_i\Big]+(\dot{s}+\dot{K}x^{(N)})dt\\
&+K\Big\{\big[(A+G){x}^{(N)}-B{R^{-1}}B^T((P+K)x^{(N)}\\&+s)+f\big]dt+\frac{1}{N}\sum_{i=1}^N\sigma dW_i\Big\}\cr
=\ &-\big[(A-\rho I)^T(Px_i+Kx^{(N)}+s)\\
&+G^T((P+K)x^{(N)}+s)\cr
&+Qx_i-\Xi x^{(N)}-\bar{\eta} \big]dt+\sum_{j=1}^N\beta_i^jdW_j.
\end{aligned}$$
This implies $\beta_i^i=\frac{1}{N}K\sigma+P\sigma$, $\beta_i^j=\frac{1}{N}K\sigma, \ j\not=i$,
\begin{align}\label{eq8a}
\rho{P}(t)=&\dot{P}(t)+A^TP(t)+P(t)A +Q \cr
-&P(t)BR^{-1}B^TP(t), \ P(T)=0,\\
\label{eq9a}
\rho K(t)\!=\!\ &\dot{K}(t)\!+\!(A+G)^TK(t)\!+\!K(t)(A+G)\!+\!G^TP(t)\cr
+&P(t)G-(P(t)+K(t))BR^{-1}B^T(P(t)+K(t))\cr
+&P(t)BR^{-1}B^TP(t)-\Xi ,\ K(T)=0,
\\ \label{eq10a}
\rho s(t)=\ &\dot{s}(t)+[A+G-BR^{-1}B^T(P+K)]^Ts(t)\cr
+&(P+K)f(t)-\bar{\eta}(t), \  s(T)=0.
\end{align}

\begin{remark}
{Note that (11) is not a standard Riccati equation. Its solvability may be referred to \cite{AFIJ03}. In particular, by Theorem 4.3 in \cite[Chapter 2]{MY99}, if $\det\Big\{[0,\ I]e^{\mathcal{A}t}\left[\begin{array}{c}
0\\
I
\end{array}\right]\Big\}>0$ with $\mathcal{A}=\!\!\left[\begin{array}{cc}
\!A-\frac{\rho}{2}I & -BR^{-1}B^T \! \\
\!-Q& -A^{T}+\frac{\rho}{2}I\!
\end{array}\right]$\!, then we have
$$P(t)=\Big\{[0,\ I]e^{\mathcal{A}t}\left[\begin{array}{c}
0\\
I
\end{array}\right]\Big\}^{-1}\Big\{[0,\ I]e^{\mathcal{A}t}
\left[\begin{array}{c}
I\\
0
\end{array}\right]\Big\}.$$
}
\end{remark}
\begin{remark}
  {Denote $\Pi=P+K$. Then from (\ref{eq8a}) and (\ref{eq9a}), $\Pi$ satisfies
\begin{equation}\label{eq11}
\begin{aligned}
&\rho \Pi(t)=\dot{\Pi}(t)+(A+G)^T\Pi(t)+\Pi(t) (A+G)\\
&-\Pi(t) BR^{-1}B^T\Pi(t)+(I-\Gamma)^TQ(I-\Gamma).
\end{aligned}
\end{equation}
with $\Pi(T)=0$. By \cite[Theorem 4.5]{SLY16}, the solvability of (\ref{eq8a}) and (\ref{eq9a}) is equivalent to the uniform convexity of two optimal control problems. Particularly, if $Q\geq 0$,
then (\ref{eq8a}) and (\ref{eq9a}) admit a unique solution, respectively.}
\end{remark}

\begin{theorem}\label{thm2}
	Assume A1) holds, and (\ref{eq8a})-(\ref{eq9a}) admit a solution, respectively. 
Then (P1) has an optimal control $$\check{u}_i(t)=-{R^{-1}}B^T[P(t)x_i(t)+K(t)x^{(N)}(t)+s(t)],$$
	where $P, K $ and $s$ are determined by (\ref{eq8a})-(\ref{eq10a}).
\end{theorem}

{To prove Theorem \ref{thm2}, we first provide a lemma, which plays a key role in the later analysis.
\begin{lemma}\label{lem2.1}
If  (\ref{eq8a}) and (\ref{eq9a}) admit a solution, respectively, then Problem (P1) is uniformly convex.
\end{lemma}
\emph{Proof.} By (\ref{eq8a}), (\ref{eq11}), and direct calculations, we have
\begin{align*}
&  \sum_{i=1}^N\mathbb{E}\int_0^Te^{-\rho t}\Big(\big\|y_i(t)-\Gamma y^{(N)}(t)\big\|^2_{Q}+\|u_i(t)\|^2_{R}\Big)dt\cr
=&\sum_{i=1}^N\mathbb{E}\int_0^Te^{-\rho t}\Big(\|y_i(t)\|^2_Q+\|y^{(N)}(t)\|^2_{\Xi }+\|u_i(t)\|^2_{R}\Big)dt\cr
=&\sum_{i=1}^N\mathbb{E}\int_0^Te^{-\rho t}\Big(\|y_i(t)-y^{(N)}(t)\|^2_Q+\|y^{(N)}(t)\|^2_{Q-\Xi }\cr
&\qquad+\|u_i(t)-u^{(N)}(t)\|^2_{R}+\|u^{(N)}(t)\|^2_{R}\Big)dt\cr
=&\sum_{i=1}^N\mathbb{E}\int_0^Te^{-\rho t}\Big(\big\|u_i(t)-u^{(N)}(t)\cr
&\qquad+R^{-1}B^TP(t)(y_i(t)-y^{(N)}(t))\big\|^2_R\cr
&\qquad+\big\|u^{(N)}(t)+R^{-1}B^T\Pi(t) y^{(N)}(t)\big\|^2_{R}\Big)dt\cr
\geq&\sum_{i=1}^N\mathbb{E}\int_0^Te^{-\rho t}\Big(\big\|u_i(t)+R^{-1}B^TP(t)(y_i(t)-y^{(N)}(t))\cr
&\qquad+R^{-1}B^T\Pi(t) y^{(N)}(t)\big\|^2_{R}\Big)dt\cr
\geq& \gamma\sum_{i=1}^N\mathbb{E}\int_0^Te^{-\rho t}\|u_i(t)\|^2dt,
\end{align*}
where the last line follows by \cite[Lemma 2.3]{SLY16}.
From Proposition \ref{thm1}, the lemma follows. \hfill{$\Box$}}

\emph{Proof of Theorem \ref{thm2}.}
Since (\ref{eq8a}) and (\ref{eq9a}) have a solution, respectively, then by \cite[Chapter 2, \S4]{MY99}, 
(\ref{eq10}) admits a unique solution. Thus, the FBSDE (\ref{eq4a}) is decoupled and the existence 
of a solution follows. From Lemma \ref{lem2.1}, (P1) is uniformly convex.
By Theorem \ref{thm1}, (P1) has an optimal control
given by
$\check{u}_i(t)=-{R^{-1}}B^T[P(t)x_i(t)+K(t)x^{(N)}(t)+s(t)],$ $t\geq 0$,
where $P, K $ and $s$ are determined by (\ref{eq8a})-(\ref{eq10a}).  \hfill{$\Box$}

As an approximation to ${x}^{(N)}$ in (\ref{eq10}), we obtain
\begin{equation}\label{eq12a}
\begin{aligned}
\frac{d\bar{x}}{dt}=&(A+G)\bar{x}(t)-B{R^{-1}}B^T(\Pi(t)\bar{x}(t)+s(t))\\
& +f(t),\ \bar{x}(0)=\bar{x}_0.
\end{aligned}
\end{equation}
Then, by Theorem \ref{thm2}, the decentralized control law for agent $i$ may be taken as
\begin{equation}\label{eq15a}
\begin{aligned}
\hat{u}_i(t)=&-{R^{-1}}B^T[P(t)\hat{x}_i(t)+K(t)\bar{x}(t)+s(t)], \\
&\ 0\leq t\leq T,
\ i=1,\cdots, N,
\end{aligned}
\end{equation}
where $P, K$, and $s$ are determined by (\ref{eq8a})-(\ref{eq10a}), and $\bar{x}$ and $\hat{x}_i$ respectively satisfy (\ref{eq12a}) and
%
%
%
\begin{align}\label{eq20}
d\hat{x}_i(t)=&\big[({A}-BR^{-1}B^TP(t))\hat{x}_i(t)+G\hat{x}^{(N)}(t)+f(t)
\cr-&BR^{-1}B^T(K(t)\bar{x}(t)+s(t))\big]dt+\sigma(t) dW_i(t).
\end{align}

\begin{remark}
Here, we firstly obtain the centralized open-loop solution by variational analysis.
  By tackling the coupled FBSDEs combined with mean field approximations, the decentralized control laws are designed. Note that in this case $s$ and $\bar{x}$ are fully decoupled and no fixed-point equation is needed.
\end{remark}


\begin{theorem}\label{thm3}
	Assume that A1) holds, and (\ref{eq8a})-(\ref{eq9a}) admit a solution, respectively. 
The set of decentralized control laws
	$\{\hat{u}_1,\cdots,\hat{u}_N\}$ in (\ref{eq15a}) has asymptotic social optimality, i.e.,
	$$\Big|\frac{1}{N}J^{\rm F}_{\rm soc}(\hat{u})-\frac{1}{N}\inf_{u\in L^2_{{\mathcal  F}}(0, T; \mathbb{R}^{nr})}J^{\rm F}_{\rm soc}(u)\Big|=O(\frac{1}{\sqrt{N}}),$$
	and the corresponding social cost is given by
  \begin{align}\label{eq16b}
J^{\rm F}_{\rm soc}(\hat{u})=&\sum_{i=1}^N\mathbb{E}\Big\{\big\|x_{i0}-x^{(N)}(0)\big\|^2_P+\big\|x^{(N)}(0)\big\|^2_{\Pi} \cr
  &+2s^T(0)x^{(N)}(0)\Big\}+Nq_T+N\epsilon_T,
  \end{align}
where
\begin{align}\label{eq17a}
  q_T=&\int_0^{T}e^{-\rho t}\big[\|\sigma (t)\|^2_{P(t)}+\| \sigma(t)\|^2_{\Pi(t)} \cr&-\|B^Ts(t)\|^2_{R^{-1}}+2s^T(t)f(t)\big]dt,\\
  \label{eq17b}
\epsilon_T=&\mathbb{E}\!\int_0^{T}\!e^{-\rho t}\|B^TK(t)(x^{(N)}(t)\!-\!\bar{x}(t))\|^2_{R^{-1}}dt.
  \end{align}
\end{theorem}

\emph{Proof.} See Appendix A. $\hfill \Box$

\section{The infinite-horizon problem}

Based on the analysis in Section \ref{sec3}, we may design the following decentralized control laws for Problem (P):
\begin{equation}\label{eq14}
\begin{aligned}
\hat{u}_i(t)=&-{R^{-1}}B^T[P\hat{x}_i(t)+(\Pi-P)\bar{x}(t)+s(t)], \\
& \ t\geq 0, \ \ i=1,\cdots, N,
\end{aligned}
\end{equation}
where $P$ and $\Pi$ are 
maximal solutions\footnote{For a Riccati equation (e.g., (\ref{eq15})), $P$ is called a maximal solution if for any solutions $P^{\prime}$, $P-P^{\prime}\geq 0$.} to the equations
\begin{align}\label{eq15}
\rho P=&A^TP+PA-PBR^{-1}B^TP+Q,\\
\label{eq16}
\rho \Pi= &(A+G)^T\Pi+\Pi (A+G)-\Pi BR^{-1}B^T\Pi+Q-\Xi ,
\end{align}
and
$s, \bar{x}\in C_{\rho/2}([0,\infty),\mathbb{R}^n)$ are determined by
\begin{align}
\rho s(t)\!=&\dot{s}(t)\!+\!(A+G\!-\!BR^{-1}B^T\Pi)^Ts(t)\!+\!\Pi f(t)-\bar{\eta}(t),\label{eq17}\\
\dot{\bar{x}}(t)=&(A+G)\bar{x}(t)-B{R^{-1}}B^T(\Pi\bar{x}(t)+s(t))\cr
&+f(t),\   \bar{x}(0)=\bar{x}_0.\label{eq18}
\end{align}
Here $s(0)$ is to be determined, and the existence conditions of $P, \Pi,s$ and $\bar{x}$ need to be investigated further.

\subsection{Uniform stabilization of subsystems}

We now list some basic assumptions for reference:

\textbf{A2)} The system $(A-\frac{\rho}{2} I, B)$ is stabilizable, and $(A+G-\frac{\rho}{2} I, B)$ is stabilizable. Particularly, $\bar{A}+G-\frac{\rho}{2}I$ is Hurwitz, where $\bar{A}\stackrel{\Delta}{=}A-BR^{-1}B^TP$.

\textbf{A3)}   $Q\geq0$, $(A-\frac{\rho}{2}I, \sqrt{Q}$) is observable, and 
$(A+G-\frac{\rho}{2}I, \sqrt{Q}(I-\Gamma))$ is observable. 

Assumptions A2) and A3) are basic in the study of the LQ optimal control problem. We will show that under some conditions, A2) is also necessary for uniform stabilization of multiagent systems.
In many cases, A3) may be weakened to the following assumption. 

{\bf{A3$^{\prime}$)}} $Q\geq0$, $(A-\frac{\rho}{2}I, \sqrt{Q}$) is detectable, and 
$(A+G-\frac{\rho}{2}I, \sqrt{Q}(I-\Gamma))$ is detectable. 


\begin{lemma}\label{lem2a}
	Under A2)-A3), (\ref{eq15}) and (\ref{eq16}) admit unique solutions $P>0, \Pi>0$, respectively, and (\ref{eq17})-(\ref{eq18}) admits a set of unique solutions $s, \bar{x}\in C_{\rho/2}([0,\infty),\mathbb{R}^n)$.
\end{lemma}

\emph{Proof.} From A2)-A3) and \cite{AM90}, (\ref{eq15}) and (\ref{eq16}) admit unique solutions $P>0, \Pi>0$ such that
$A-BR^{-1}B^TP-\frac{\rho}{2}I$ and $A+G-BR^{-1}B^T\Pi-\frac{\rho}{2}I$ are Hurwitz, respectively. From an argument in \cite[Appendix A]{WZ12}, we obtain
$s\in C_{\rho/2}([0,\infty),\mathbb{R}^n)$ if and only if 
$$s(t)=\int_t^{\infty}e^{-(A+G-BR^{-1}B^T\Pi-{\rho}I)(t-\tau)}(\Pi f(\tau)-\bar{\eta}(\tau))d\tau.$$
$\hfill \Box$


%

\begin{lemma}\label{lem2}
	Let A1)-A3) hold. Then for Problem (P),
	\begin{equation}\label{eq13a}
	\mathbb{E}\int_0^{\infty} e^{-\rho t}\|\hat{x}^{(N)}(t)-\bar{x}(t)\|^2dt=O(\frac{1}{N}),
	\end{equation}
	
	where $\hat{x}^{(N)}=\sum_{i=1}^N\hat{x}_i$, and $\bar{x}$ satisfies (\ref{eq18}).
\end{lemma}

\emph{Proof.} See Appendix \ref{app b}.  $\hfill \Box$

It is shown that the decentralized control laws (\ref{eq15a}) uniformly stabilize the systems (\ref{eq1}) .

\begin{theorem}\label{thm4}
	Let A1)-A3) hold.  Then for any $N$, 
	\begin{equation}\label{eq13b}
	\sum_{i=1}^N\mathbb{E}\int_0^{\infty} e^{-\rho t} \left(\|\hat{x}_i(t)\|^2+\|\hat{u}_i(t)\|^2\right)dt<\infty.
	\end{equation}

\end{theorem}

\emph{Proof.} See Appendix \ref{app b}.  $\hfill \Box$

We now give two equivalent conditions for uniform stabilization of multiagent systems.

\begin{theorem}\label{thm5}
	Let A3) hold. Assume that (\ref{eq15})-(\ref{eq16}) admit symmetric solutions. Then for Problem (P) the following statements are equivalent:
	
	(i)  For any initial condition $(\hat{x}_1(0),\cdots, \hat{x}_N(0))$ satisfying A1),
	\begin{equation}\label{eq23}
	\sum_{i=1}^N\mathbb{E}\int_0^{\infty} e^{-\rho t} \left(\|\hat{x}_i(t)\|^2+\|\hat{u}_i(t)\|^2\right)dt<\infty.
	\end{equation}
	
	(ii) Equations (\ref{eq15}) and (\ref{eq16}) admit unique maximal solutions 
such that $P>0, \Pi>0$, and
 $\bar{A}+G-\frac{\rho}{2}I$ is Hurwitz.
	
	(iii) A2) holds.
	
\end{theorem}

\emph{Proof.} See the Appendix C.  $\hfill \Box$


For $G=0$, we have a simplified version of Theorem \ref{thm5}.
\begin{corollary}
	Assume that A3) holds and $G=0$. Assume that (\ref{eq15})-(\ref{eq16}) admit  symmetric solutions. Then the following statements are equivalent:
	
	(i)  For any $(\hat{x}_1(0),\cdots, \hat{x}_N(0))$ satisfying A1),
	\begin{equation*}
	\sum_{i=1}^N\mathbb{E}\int_0^{\infty} e^{-\rho t} \left(\|\hat{x}_i(t)\|^2+\|\hat{u}_i(t)\|^2\right)dt<\infty.
	\end{equation*}	
	
	(ii) Equations (\ref{eq15}) and (\ref{eq16}) admit unique maximal solutions such that $P>0, \Pi>0$, respectively.
	
	(iii) A2) holds.
\end{corollary}

When A3) is weakened to A3$^{\prime}$), we have the following equivalent conditions of uniform stabilization. 
\begin{theorem}\label{thm6}
	Let A3$^{\prime}$) hold. Assume that (\ref{eq15})-(\ref{eq16}) admit solutions. Then the following are equivalent:
	
	(i)  For any initials $(\hat{x}_1(0),\cdots, \hat{x}_N(0))$ satisfying A1),
	\begin{equation*}
	\sum_{i=1}^N\mathbb{E}\int_0^{\infty} e^{-\rho t} \left(\|\hat{x}_i(t)\|^2+\|\hat{u}_i(t)\|^2\right)dt<\infty.
	\end{equation*}
	
	(ii) Equations (\ref{eq15}) and (\ref{eq16}) admit unique maximal solutions $P\geq0, \Pi\geq0$, and $\bar{A}+G-\frac{\rho}{2}I$ is Hurwitz.
	
	(iii) A2) holds.
	
\end{theorem}

\emph{Proof.} See the Appendix C.  $\hfill \Box$


For the more general case that $Q$ are indefinite, we have the following equivalent
conditions for uniform stabilization of all the subsystems.
Assume

$\mathbf{A3^{\prime\prime})}$ both $M_1$ and $M_2$ have no eigenvalues on the imaginary axis, where
$$M_1=\left[\begin{array}{cc}
A-\frac{\rho}{2}I & BR^{-1}B^T \\
Q& -A^{T}+\frac{\rho}{2}I
\end{array}\right],$$
$$M_2=\left[\begin{array}{cc}
A+G-\frac{\rho}{2}I & BR^{-1}B^T \\
Q-\Xi & -(A+G)^{T}+\frac{\rho}{2}I
\end{array}\right].$$

\begin{theorem}\label{thm7}
	Assume that ${A3^{\prime\prime})}$ holds, and (\ref{eq15})-(\ref{eq16}) admit solutions. Then the following  are equivalent:
	
	(i)  For any $(\hat{x}_1(0),\cdots, \hat{x}_N(0))$ satisfying A1),
	\begin{equation*}
	\sum_{i=1}^N\mathbb{E}\int_0^{\infty} e^{-\rho t} \left(\|\hat{x}_i(t)\|^2+\|\hat{u}_i(t)\|^2\right)dt<\infty.
	\end{equation*}	
	
	(ii) Equations (\ref{eq15}) and (\ref{eq16}) admit unique $\rho$-stabilizing solutions\footnote{For a Riccati equation (\ref{eq15}), $P$ is called a $\rho$-stabilizing solution if $P$
		satisfies (\ref{eq15}) and all the eigenvalues of $A-BR^{-1}B^TP-\frac{\rho}{2}I$ are in left half-plane.} (which are also the maximal solutions), and $\bar{A}+G-\frac{\rho}{2}I$ is Hurwitz.
	
	(iii) A2) holds.
\end{theorem}

\begin{remark}
$M_1$ and $M_2$ are Hamiltonian matrices. The Hamiltonian matrix plays a significant role in studying general algebraic Riccati equations. See more details of the property of Hamiltonian matrices in \cite{AFIJ03}, \cite{M77}.
\end{remark}
\begin{remark}
 For the case $Q=0$ and $G=0$, the Hamiltonian matrices reduce to
$$M_1=M_2=\left[\begin{array}{cc}
A-\frac{\rho}{2}I & BR^{-1}B^T \\
0& -A^{T}+\frac{\rho}{2}I
\end{array}\right].$$
Then it follows from Theorem 4.4 that if $A-\frac{\rho}{2}I$ have no eigenvalues on the imaginary axis, the decentralized controls (\ref{eq15a}) uniformly stabilize the systems (\ref{eq1}) if and only if $(A-\frac{\rho}{2}I, B)$ is stabilizable. Since $Q=0$ and $A-\frac{\rho}{2}I$ is not Hurwitz necessarily, the system ($A-\frac{\rho}{2}I,\sqrt{Q}$) is not detectable, which implies that the assumptions of Theorem 4.3 in [21] does not hold.
\end{remark}

To show Theorem \ref{thm7}, we need two lemmas. The first lemma is copied from \cite[Theorem 6]{M77}.
\begin{lemma}\label{lem5}
	Equations (\ref{eq15}) and (\ref{eq16}) admit unique $\rho$-stabilizing solutions (which are also the maximal solutions) if and only if
	A2) and ${A3^{\prime\prime})}$ hold.  
\end{lemma}

\begin{lemma}\label{lem4}
	Let A1) hold. Assume that (\ref{eq15}) and (\ref{eq16}) admit $\rho$-stabilizing solutions, respectively, and $\bar{A}+G-\frac{\rho}{2}I$ is Hurwitz. Then
	\begin{equation*}
	\sum_{i=1}^N\mathbb{E}\int_0^{\infty} e^{-\rho t} \left(\|\hat{x}_i(t)\|^2+\|\hat{u}_i(t)\|^2\right)dt<\infty.
	\end{equation*}
\end{lemma}
\emph{Proof.} From the definition of $\rho$-stabilizing solutions, $A-BR^{-1}B^TP-\frac{\rho}{2}I$ and $A+G-BR^{-1}B^T\Pi-\frac{\rho}{2}I$ are Hurwitz. By the argument in the proof of Theorem \ref{thm4}, the lemma follows. $\hfill \Box$

\emph{Proof of Theorem \ref{thm7}.}  By using Lemmas \ref{lem5} and \ref{lem4} together with a similar argument in the proof of Theorem \ref{thm4}, 
the theorem follows. $\hfill \Box$

\begin{example}\label{ex1}
	Consider a scalar system with $A=a$, $B=b$, $G=g$, $Q=q$, $\Gamma=\gamma$, $R=r>0$. Then
	$$M_1=\left[\begin{array}{cc}
	a-{\rho}/{2} & b^2/r \\
	q& -a+{\rho}/{2}
	\end{array}\right],$$
	$$M_2=\left[\begin{array}{cc}
	a+g-{\rho}/{2} & {b^2}/{r} \\
	q(1-\gamma)^2& -(a+g-{\rho}/{2})
	\end{array}\right].$$
	By direct computations, neither $M_1$ nor $M_2$ has eigenvalues in the imaginary axis if and only if
	\begin{align}\label{eq41a}
	&(a-\frac{\rho}{2})^2+\frac{b^2}{r}q>0,\\ \label{eq42a}
	&(a+g-\frac{\rho}{2})^2+\frac{b^2}{r}(1-\gamma)^2q>0.
	\end{align}
	Note that if $q>0$ (or $a-{\rho}/{2}<0$,\ $q=0$), i.e., $(a-{\rho}/{2},\sqrt{q})$ is observable (detectable), then (\ref{eq41a}) holds,
	and  if $(1-\gamma)^2q>0$ ($a+g-{\rho}/{2}<0,$\ $q=0$), i.e., $(a+g-{\rho}/{2},\sqrt{q}(1-\gamma))$ is observable (detectable), then (\ref{eq42a}) holds.
	
	For this model, the Riccati equation (\ref{eq15}) is written as
	\begin{equation}\label{eq43a}
	\frac{b^2}{r}p^2-(2a-\rho)p-q=0.
	\end{equation}
	Let $\Delta=4[(a-{\rho}/{2})^2+{b^2q}/{r}]$. If (\ref{eq41a}) holds then $\Delta>0$, which implies (\ref{eq43a}) admits two solutions. If $q>0$ then (\ref{eq43a}) has a unique positive solution such that $a-b^2p/r-{\rho}/{2}=-\sqrt{\Delta}/2<0$.
	If $q=0$ and $a-\rho/2<0$ then (\ref{eq43a}) has a unique non-negative solution $p=0$ such that $a-b^2p/r-{\rho}/{2}=a-{\rho}/2<0$.
	
	Assume that (\ref{eq41a}) and (\ref{eq42a}) hold. By Theorem \ref{thm7}, the system is uniformly stable if and only if
	$(a-\rho/2,b)$ is stabilizable (i.e., $b\not=0$ or $a-\rho/2<0$), and $a-b^2p/r-{\rho}/{2}+g<0$. Note that $a-b^2p/r-{\rho}/{2}<0$. When $g\leq0$,
we have $a-b^2p/r-{\rho}/{2}+g<0$.
	
\end{example}

\begin{example}
	We further consider the model in Example \ref{ex1} for the case that $a+g=\rho/2$ and $\gamma=1$ (i.e., (\ref{eq42a}) does not hold).  In this case, the Riccati equation (\ref{eq16}) admits a unique solution $\Pi=0$.  (\ref{eq17}) becomes
	$\rho s(t)=\dot{s}(t)+\frac{\rho}{2}s(t)$
	and has a unique solution $s(t)\equiv0$ in $C_{\rho/2}([0,\infty),\mathbb{R})$. Thus, $\bar{x}$ satisfies
	\begin{equation}\label{eq44}
	\frac{d\bar{x}}{dt}=\frac{\rho}{2}\bar{x}(t)+f(t).
	\end{equation}
	Assume that $f$ is a constant.
	Then (\ref{eq44}) does not admit a solution in $C_{\rho/2}([0,\infty),\mathbb{R})$ unless $\bar{x}(0)=-{2f}/{\rho}$.
	
\end{example}

\subsection{Asymptotic social optimality}
Now we are in a position to state the asymptotic optimality of
the decentralized control. 
\begin{theorem}\label{thm8}
	Let A1)-A3) hold. For Problem (P), the set of decentralized control laws
	$\{\hat{u}_1,\cdots,\hat{u}_N\}$ given by (\ref{eq14}) has asymptotic social optimality, i.e.,
	$$\Big|\frac{1}{N}J_{\rm soc}(\hat{u})-\frac{1}{N}\inf_{u\in \mathcal{U}_c}J_{\rm soc}(u)\Big|=O({1}/{\sqrt{N}}).
$$
\end{theorem}

\emph{Proof.} We first prove that for $u\in \mathcal{U}_c$, $J_{\rm soc}(u)< NC_1$ implies that
\begin{equation}\label{eq36a}
\sum_{i=1}^N\mathbb{E}\int_0^{\infty}e^{-\rho t}(\|x_i(t)\|^2+\|u_i(t)\|^2)dt<NC_2,
\end{equation} for all $i=1,\cdots,N$.
From $J_{\rm soc}(u)< NC_1$, we have
$\sum_{i=1}^N\mathbb{E}\int_0^{\infty}e^{-\rho t}\|u_i(t)\|^2dt<NC$ and
\begin{equation}\label{eq36}
\sum_{i=1}^N\mathbb{E}\int_0^{\infty}e^{-\rho t}\big\|x_i(t)-\Gamma x^{(N)}(t)\big\|^2_Qdt<NC,
\end{equation}
which further implies that
\begin{equation}\label{eq37}
\begin{aligned}
&\mathbb{E}\int_0^{\infty}e^{-\rho t}\big\|(I-\Gamma)x^{(N)}(t)\big\|^2_Q \\
\leq&\frac{1}{N}\sum_{i=1}^N
\mathbb{E}\int_0^{\infty}e^{-\rho t}\big\|x_i(t)-\Gamma x^{(N)}(t)\big\|^2_Q dt<C.
\end{aligned}
\end{equation}
By (\ref{eq1}) we have
\begin{equation*}
\begin{aligned}
dx^{(N)}(t)=\ & \left[(A+G)x^{(N)}(t)+Bu^{(N)}(t)+f(t)\right]dt\\
&+\frac{1}{N}\sum_{i=1}^N\sigma(t) dW_i(t),
\end{aligned}
\end{equation*}
which leads to for any $r\in [0,1]$,
\begin{equation}\label{eq39}
\begin{aligned}
x^{(N)}(t)=\ &e^{(A+G)r}x^{(N)}(t-r)\\
&+\int_{t-r}^te^{(A+G)(t-\tau)}[Bu^{(N)}(\tau)+f(\tau)]d\tau\\
&+\frac{1}{N}\sum_{i=1}^N\int_{t-r}^te^{(A+G)(t-\tau)}\sigma(\tau) dW_i(\tau).
\end{aligned}
\end{equation}
By $J_{\rm soc}(u)< C_1$ and basic SDE estimates, we can find a constant $C$ such that
$$\mathbb{E}\int_r^{\infty}e^{-\rho t}\Big\|\int_{t-r}^te^{(A+G)(t-\tau)}Bu^{(N)}(\tau)d\tau\Big\|^2dt\leq C.$$
From (\ref{eq37}) and (\ref{eq39}) we obtain
$$
\begin{aligned}
\mathbb{E}\int_r^{\infty}&e^{-\rho t}[x^{(N)}(t-r)]^Te^{(A+G)^Tr}(I-\Gamma)^TQ(I-\Gamma)\\
&\cdot e^{(A+G)r}x^{(N)}(t-r)dt\leq C,
\end{aligned}
$$
which implies that for any $r\in [0,1]$,
$$
\begin{aligned}
\mathbb{E}\int_0^{\infty}&e^{-\rho t}[x^{(N)}(\tau)]^Te^{(A+G)^Tr}(I-\Gamma)^TQ(I-\Gamma)\\
&\cdot e^{(A+G)r}x^{(N)}(\tau)d\tau\leq C.
\end{aligned}
$$
By taking integration with respect to $r\in [0,1]$, we obtain
$$
\begin{aligned}
\mathbb{E}\int_0^{\infty}&e^{-\rho t}[x^{(N)}(\tau)]^T\Big[\int_0^1e^{(A+G)^Tr}(I-\Gamma)^TQ(I-\Gamma)\\
&\cdot e^{(A+G)r}dr\Big]x^{(N)}(\tau)d\tau\leq C.
\end{aligned}
$$
This together with A3) lead to
 \begin{equation}\label{eq41}
\mathbb{E}\int_0^{\infty}e^{-\rho t}\|x^{(N)}(t)\|^2dt<C,
\end{equation}
which with (\ref{eq36}) further gives
\begin{equation}\label{eq41c}
\sum_{i=1}^N\mathbb{E}\int_0^{\infty}e^{-\rho t}\|x_i(t)\|^2_Qdt<NC.
\end{equation}

By (\ref{eq1}), we have
\begin{equation}\label{eq42}
\begin{aligned}
x_i(t)=\ &e^{Ar}x_i(t-r)\\
&+\int_{t-r}^te^{A(t-\tau)}[Bu_i(\tau)+f(\tau)+Gx^{(N)}(\tau)]d\tau\\
&+\int_{t-r}^te^{A(t-\tau)}\sigma(\tau) dW_i(\tau).
\end{aligned}
\end{equation}
It follows from (\ref{eq41}) that
\begin{align*}
&\mathbb{E}\int_r^{\infty}e^{-\rho t}\Big\|\int_{t-r}^te^{A(t-\tau)}Gx^{(N)}(\tau)d\tau\Big\|^2dt\cr
\leq & \mathbb{E}\int_0^{\infty}e^{-\rho \tau}\|Gx^{(N)}(\tau)\|^2\int_{0}^r\big\|e^{(A-\frac{\rho}{2}I)v}\big\|^2dvd\tau\leq C.
\end{align*}
From (\ref{eq41c}) and (\ref{eq42}), we obtain that
$$\sum_{i=1}^N\mathbb{E}\int_r^{\infty}e^{-\rho t}x_i^T(t-r)e^{A^Tr}Qe^{Ar}x_i(t-r)dt\leq NC.$$
This together with A3)
implies that
\begin{equation*}
\sum_{i=1}^N\mathbb{E}\int_0^{\infty}e^{-\rho t}\|x_i(t)\|^2dt<NC,
\end{equation*}
which gives (\ref{eq36a}). From this with Theorem \ref{thm4},
\begin{equation*}
\sum_{i=1}^N\mathbb{E}\int_0^{\infty}e^{-\rho t}\big(\|\tilde{x}_i(t)\|^2+\|\tilde{u}_i(t)\|^2\big)dt<NC.
\end{equation*}
By a similar argument to the proof of Theorem \ref{thm3} combined with Lemma \ref{lem2}, the conclusion follows.
$\hfill \Box$

If A3) is replaced by A3$^{\prime}$), the decentralized control (\ref{eq14}) still has asymptotic social optimality.
\begin{corollary}\label{cor4.2}
	Assume that A1)-A2), A3$^{\prime}$) hold.
	The decentralized control (\ref{eq14}) is asymptotically social optimal.
\end{corollary}
\emph{Proof.}   Without loss of generality, we simply assume $A+G=\hbox{diag}\{\mathbb{A}_{1},\mathbb{A}_2\}$, where $\mathbb{A}_1-(\rho/2) I$ is Hurwitz, and $-(\mathbb{A}_2-(\rho/2) I)$ is Hurwitz (If necessary, we may apply a nonsingular linear transformation as in the proof of Theorem \ref{thm6}). Write $x^{(N)}=[z_1^T,z_2^T]$ and $Q^{1/2}(I-\Gamma)=[S_1,S_2]$ such that
$\big\|(I-\Gamma)x^{(N)}(t)\big\|_Q^2=
\|S_1{z}_1(t)+S_2z_2(t)\|^2,$
and $(\mathbb{A}_2-(\rho/2) I,S_2)$ is observable. 
By the proof of Theorem \ref{thm4} or \cite{H10}, $\mathbb{E}\int_0^{\infty}e^{-\rho t}\|u^{(N)}(t)\|^2dt<\infty$ implies $\mathbb{E}\int_0^{\infty}e^{-\rho t}\|z_1(t)\|^2dt<\infty$, which together with (\ref{eq37}) gives  $\mathbb{E}\int_0^{\infty}e^{-\rho t}\|S_2z_2(t)\|^2dt<\infty$. This and the observability of $(A_2-(\rho/2) I,S_2)$ leads to $\mathbb{E}\int_0^{\infty}e^{-\rho t}\|z_2(t)\|^2dt<\infty$. Thus, $\mathbb{E}\int_0^{\infty}e^{-\rho t}\|x^{(N)}(t)\|^2dt$ \\$<\infty$.  The other parts of the proof are similar to that of Theorem \ref{thm8}.
$\hfill \Box$

For the case that $Q$ are indefinite, we have the following result of asymptotic optimality.

\begin{theorem}\label{thm4.6}
	Let A1)-A2), ${A3^{\prime\prime})}$ hold. Assume (\ref{eq15})-(\ref{eq16}) admit negative definite solutions
$P^{-}<0$ and $\Pi^-<0$, respectively.
	Then, the set of decentralized control in (\ref{eq14}) is asymptotically socially optimal. Furthermore, if $\{x_{i0}\}$ have the same variance, then the  asymptotic average social optimum is given by
$$\lim_{N\to\infty}\!\frac{1}{N}J_{\rm soc}(\hat{u})\!=\!\mathbb{E}\big[\|x_{i0}-\bar{x}_0\|_P^2+\|\bar{x}_0\|^2_{\Pi}+2s^T(0)\bar{x}_0\big]\!+q_{\infty},$$
where
\begin{align}
\label{eq39a}
  q_{\infty}=&\int_0^{\infty}e^{-\rho t}\big[\|\sigma (t)\|^2_{P}+\| \sigma(t)\|^2_{\Pi} \cr&-\|B^Ts(t)\|^2_{R^{-1}}+2s^T(t)f(t)\big]dt.
\end{align}
\end{theorem}
\emph{Proof.} From the above assumptions and Theorem \ref{thm7}, the Riccati equation (\ref{eq15}) admits a $\rho$-stabilizing solution $P$ and a negative definite solution $P^-$;
(\ref{eq16}) has a $\rho$-stabilizing solution $\Pi$ and a negative definite solution $\Pi^-$.
By a similar argument in the proof of Lemma \ref{lem2.1}, we obtain for any $u\in {\mathcal  U}_c$,
\begin{align*}
&J_{\rm soc}(u)\cr
  =&\sum_{i=1}^N\mathbb{E}\int_0^{\infty}e^{-\rho t}\Big(\|{x}_i-{x}^{(N)}\|^2_Q+\|{x}^{(N)}\|^2_{Q-\Xi }+\|\eta\|_Q^2\cr
&-2\eta^TQ(I-\Gamma){x}_i+\|{u}_i-{u}^{(N)}\|^2_{R}+\|{u}^{(N)}\|^2_{R}\Big)dt\cr
=&\sum_{i=1}^N\mathbb{E}\big[\|x_{i0}-{x}^{(N)}(0)\|^2_{P^{-}}+\|{x}^{(N)}(0)\|^2_{\Pi^-}\cr&
+2s^T(0){x}^{(N)}(0)\big]
-\lim_{T\to \infty}\sum_{i=1}^Ne^{-\rho T}\mathbb{E}\big[\|{x}^{(N)}(T)\|^2_{\Pi^{-}}\cr
&+\|x_{i}(T)-{x}^{(N)}(T)\|^2_{P^{-}}+2s^T(T){x}^{(N)}(T)\big]\cr
&+\sum_{i=1}^N\mathbb{E}\int_0^{\infty}e^{-\rho t}\Big(\big\|{u}^{(N)}+R^{-1}B^T\Pi^-{x}^{(N)}\big\|^2_{R}\cr
&+\big\|{u}_i-{u}^{(N)}+R^{-1}B^TP^-({x}_i-{x}^{(N)})\big\|^2_R\Big)dt+  q_{\infty}.
\end{align*}
By \cite[Theorem 8]{W71}, the centralized optimal control exists and the optimal state is $\rho$-stable.
Thus, we only need to consider the following control set
$$\begin{aligned}
  \mathcal{U}_c^{\prime}=&\Big\{(u_1,\cdots,u_N)|u_i(t) \hbox{ is adapted to } \mathcal{F}_t,\\
&  \qquad \mathbb{E}\int_0^{\infty}e^{-\rho t}\|x_i(t)\|^2dt<\infty, \forall i \Big\}.
  \end{aligned}$$
For any $u\in {\mathcal  U}_c^{\prime}$ satisfying $J_{\rm soc}(u)\leq NC$, we have
\begin{align}\label{eq38b}
&J_{\rm soc}(u)\cr
=&\sum_{i=1}^N\mathbb{E}\big[\|x_{i0}-{x}^{(N)}(0)\|^2_{P}+\|{x}^{(N)}(0)\|^2_{\Pi}+2s^T(0)\bar{x}_0\big]
\cr
&+\!\sum_{i=1}^N\mathbb{E}\int_0^{\infty}\!\!e^{-\rho t}\Big(\big\|{u}_i\!-\!{u}^{(N)}\!+\!R^{-1}B^TP({x}_i\!-\!{x}^{(N)})\big\|^2_R\cr
&+\big\|{u}^{(N)}+R^{-1}B^T\Pi{x}^{(N)}\big\|^2_{R}\Big)dt+  q_{\infty}\leq NC.\cr
\end{align}
Denote $v^{(N)}={u}^{(N)}+R^{-1}B^T\Pi{x}^{(N)}$. From (\ref{eq1}),
$$\begin{aligned}
  dx^{(N)}(t)=&(A+G-BR^{-1}B^T\Pi)x^{(N)}(t)dt\cr
  &+Bv^{(N)}(t)dt+\frac{1}{N}\sum_{i=1}^N\sigma(t) dW_i(t).
\end{aligned}$$
By \cite{H10},  there exists $C_1,C_2>0$ such that
$$\mathbb{E}\int_0^{\infty}e^{-\rho t}\|x^{(N)}\|^2dt\leq C_1 \mathbb{E}\int_0^{\infty}e^{-\rho t}\|v^{(N)}\|^2+C_2.$$
This together with (\ref{eq38b}) gives
\begin{align}\label{eq39b}
  &\sum_{i=1}^N\mathbb{E}\int_0^{\infty}e^{-\rho t}(\|x^{(N)}\|^2+\|u^{(N)}\|^2)dt
 \cr=&N\mathbb{E}\int_0^{\infty}\!e^{-\rho t}(\|x^{(N)}\|^2+\|v^{(N)}-R^{-1}B^T\Pi{x}^{(N)}\|^2)dt\cr
  \leq& NC_3\mathbb{E}\int_0^{\infty}e^{-\rho t}\|v^{(N)}\|^2+NC_4\leq NC.
\end{align}
Similarly, we have
$$\sum_{i=1}^N\mathbb{E}\int_0^{\infty}e^{-\rho t}(\|x_i-x^{(N)}\|^2+\|u_i-u^{(N)}\|^2)dt\leq NC.$$
From this and (\ref{eq39b}),
$$\sum_{i=1}^N\mathbb{E}\int_0^{\infty}e^{-\rho t}(\|x_i\|^2+\|u_i\|^2)dt\leq NC.$$
The remainder of the proof can follow by that of Theorem \ref{thm3}.  For the case that $\{x_{i0}\}$ have the same variance, from (\ref{eq16b}), the asymptotic average social optimum ($\lim_{N\to\infty}\frac{1}{N}J_{\rm soc}(\hat{u})$) is given by
$\mathbb{E}\big[\|x_{i0}-\bar{x}_0\|_P^2+\|\bar{x}_0\|^2_{\Pi}+2s^T(0)\bar{x}_0\big]+q_{\infty}.$
$\hfill \Box$

\begin{remark}
  The work \cite{HCM12} investigated mean field LQ problem (P) with $Q\geq 0$. To obtain asymptotic social optimality, they need
  $Q>0$ and $I-\Gamma$ is nonsingular. In Corollary \ref{cor4.2}, we have loosed the assumption to A3$^{\prime})$, i.e., $(A- (\rho/2) I, \sqrt{Q})$
  and  $(A- (\rho/2) I, \sqrt{Q}(I-\Gamma))$ are detectable.  In Theorem \ref{thm4.6}, we further give the condition for the case of indefinite $Q$. Particularly,
  for the scalar case, the condition is equivalent to (\ref{eq41a})-(\ref{eq42a}).  It can be verified that the assumption $Q>0$ and $I-\Gamma$ is nonsingular implies
   (\ref{eq41a})-(\ref{eq42a}), but the converse is not true.
\end{remark}


\subsection{Comparison to previous solutions}

In this section, we compare the proposed decentralized control laws 
with the feedback decentralized strategies in previous works. 

We first introduce a definition from \cite{BO82}.
\begin{definition}
	For a control problem with an admissible control set $\mathcal{U}$, a control law $u\in \mathcal{U}$ is said to be a representation of another control
	$u^*\in \mathcal{U}$ if
	
	(i) they both generate the same unique state trajectory, and
	
	(ii) they both have the same open-loop value on this trajectory.
\end{definition}

For Problem (P), let $f=0$, and $G=0$.
In \cite[Theorem 4.3]{HCM12}, the decentralized control laws are given by
\begin{equation}\label{eq47}
\breve{u}_i(t)=-R^{-1}B^T(Px_i(t)+\bar{s}(t)),\quad i=1,\cdots,N,
\end{equation}
where $P$ is the stabilizing solution of (\ref{eq15}), and $\bar{s}=\bar{K}{x}^{\dag}+\phi.$ Here $\bar{K}$ satisfies
\begin{align*}
  \rho\bar{K}=\ &\bar{K}\bar{A}+\bar{A}^T\bar{K}
-\bar{K}BR^{-1}B^T\bar{K}^T-\Xi ,
\end{align*}
and ${x}^{\dag}, \phi\in C_{\rho/2}([0,\infty),\mathbb{R}^n) $ are given by
\begin{align*}
\frac{{d\bar{x}}^\dag}{dt}=\ &\bar{A}\bar{x}^{\dag}(t)-BR^{-1}B^T(\bar{K}\bar{x}^\dag(t)+\phi(t)),\bar{x}^\dag(0)=\bar{x}_0,\cr
\frac{d\phi}{dt}=\ &-[A-BR^{-1}B^T(P+\bar{K})-{\rho}I]\phi(t)+\bar{\eta}(t),
\end{align*}
in which $\bar{A}=A-BR^{-1}B^TP$ and $\phi(0)$ is to be determined by $\phi\in C_{\rho/2}([0,\infty),\mathbb{R}^n)$.
By comparing this with (\ref{eq16})-(\ref{eq18}), one can obtain that $\bar{K}=\Pi-P$, $\bar{x}=\bar{x}^{\dag}$ and $\phi=s$. From the above discussion, we have the equivalence of the two sets of decentralized control laws.
\begin{proposition}
	The set of decentralized control laws $\{\hat{u}_1,\cdots,\hat{u}_N\}$ in (\ref{eq14}) is
	a representation of $\{\breve{u}_1,\cdots,\breve{u}_N\}$ given by (\ref{eq47}).
\end{proposition}

\begin{remark}
The work \cite{HCM12} studied the problem (P) with $Q\geq0$ by the fixed-point approach.
 In Theorem 4.3, they have shown that the fixed-point equation admits a unique solution, when $(A- (\rho/2) I,\sqrt{Q})$ is detectable and
 $\Xi=\Gamma^TQ +Q\Gamma-\Gamma^TQ\Gamma \leq 0$. In fact, the above assumption 
 is merely a sufficient condition to ensure $A3^{\prime})$ ($A- (\rho/2) I,\sqrt{{Q}-\Xi)}$ is detectable). 
\end{remark}

\begin{remark}
The work \cite{HZ19} investigated asymptotic solvability of mean field LQ games by the re-scaling method. They considered (\ref{eq1})-(\ref{eq2}) with $Q\geq0$ and derived a low-dimensional ordinary differential equation system
by dynamic programming.
Actually, the method proposed in this paper can be viewed as a type of direct approach. Different from \cite{HZ19}, we tackle directly high-dimensional FBSDEs along the line of maximum principle.
\end{remark}

\section{Numerical Examples}
Now, two numerical examples are given to illustrate the effectiveness of the proposed decentralized control.

\begin{figure}[H]
	\centering
	\includegraphics[width=0.8\linewidth]{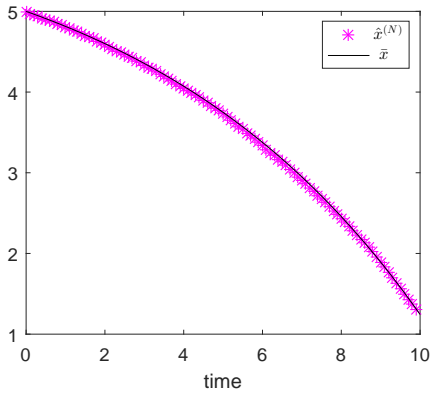}
	\caption{Curves of $\bar{x}$ and $\hat{x}^{(N)}$.}
	\label{consistency2}
\end{figure}

We first consider a scalar system with $30$ agents in Problem (P). Take $A=0.8, B=R= 1,{Q=-0.1}, G=-0.2,f(t)=1,\sigma(t)=0.2, \rho=0.6, \Gamma=0.2,\eta=5 $ in (\ref{eq1})-(\ref{eq2}). The initial states of $50$ agents are taken independently from a normal distribution $N(5,0.3)$.  Note that $B\not=0$, and $\bar{A}+G-\frac{\rho}{2}I=  -0.5873<0$. Then
A1)-A2) hold. Since $M_1=\left[\begin{array}{cc}
0.5 & 1 \\
-0.1& -0.5
\end{array}\right],
M_2=\left[\begin{array}{cc}
0.3 & 1 \\
-0.064& -0.3
\end{array}\right]$ have no eigenvalues on the imaginary axis, A3$^{\prime\prime})$ also holds. Under the control law (\ref{eq14}), 
the trajectories of $\bar{x}$ and $\hat{x}^{(N)}$ in Problem (P) are shown in Fig. \ref{consistency2}.
It can be seen that $\bar{x}$ and $\hat{x}^{(N)}$ coincide well, which illustrate the consistency of mean field approximations.


\begin{figure}[H]
	\centering
	\includegraphics[width=0.8\linewidth]{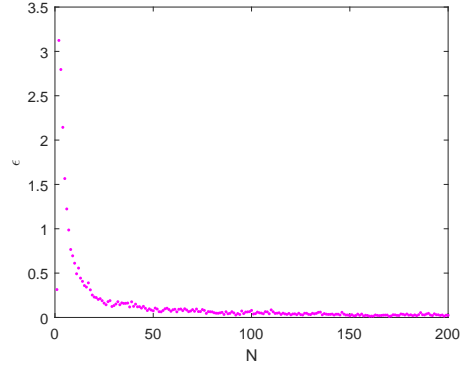}
	\caption{Curves of $\epsilon$ with resect to $N$.}
	\label{epsilon2}
\end{figure}

Denote $\epsilon=\Big|\frac{1}{N}J_{\rm soc}(\hat{u})-\frac{1}{N}\inf_{u\in \mathcal{U}_c}J_{\rm soc}(u)\Big|$. By Theorems \ref{thm3} and \ref{thm4.6}, we obtain
$\epsilon=\int_0^{\infty}e^{-\rho t}\|B^TK(x^{(N)}(t)-\bar{x}(t))\|^2_{R^{-1}}dt.$
The cost gap $\epsilon$ is demonstrated in Fig. \ref{epsilon2} where the agent number $N$ grows from 1 to 200.

\begin{figure}[H]
	\centering
	\includegraphics[width=0.85\linewidth]{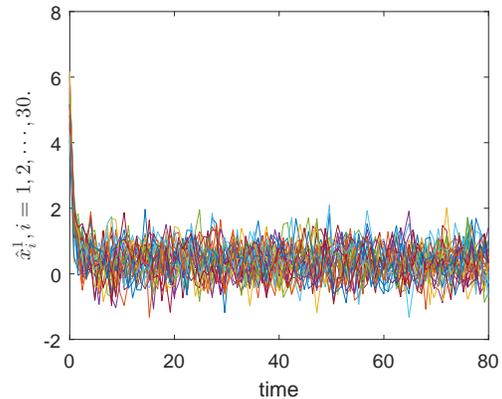}
	\caption{Curves of $\hat{x}^1_i$, $i=1,\cdots,30$.}
	\label{vector1}
\end{figure}

\begin{figure}[H]
	\centering
	\includegraphics[width=0.9\linewidth]{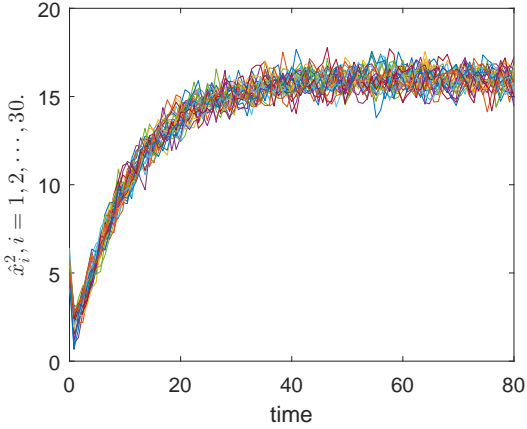}
	\caption{Curves of $\hat{x}^2_i$, $i=1,\cdots,30$.}
	\label{vector2}
\end{figure}

Finally, we consider the 2-dimensional case of Problem (P). Take parameters as follows:
$A = \left[\begin{array}{cc}
0.1 & 0\\
-1  & 0.2\\
\end{array}\right]$,
$B = \left[\begin{array}{c}
1 \\
1\\
\end{array}\right]$,
$G = \left[\begin{array}{cc}
-0.5 & 0\\
0  & -0.3\\
\end{array}\right]$,
$Q = \left[\begin{array}{cc}
1\ &\ 0\\
0\  &\ 1\\
\end{array}\right]$,
$\Gamma = \left[\begin{array}{cc}
1\ &\ 0\\
1\  &\ 1\\
\end{array}\right]$,
$R = \left[\begin{array}{cc}
1\ &\ 0\\
0\  &\ 1\\
\end{array}\right]$,
$\eta = \left[\begin{array}{c}
0\\
0.5\\
\end{array}\right]$, $f = [1\ \ 1]^T$ and $\sigma = [0.5\ \ 0.5]^T$. Denote
$\hat{x}_i(t)=[\hat{x}^1_i(t)\ \hat{x}^2_i(t)]^T$.
 Both of $\hat{x}^1_i(0)$ and $\hat{x}^2_i(0)$ are taken independently from a normal distribution $N(5,0.5)$. Under the control laws (\ref{eq14}), the trajectories of $\hat{x}^1_i$ and $\hat{x}^2_i$, $i=1,\cdots,N$ are shown in Figs. \ref{vector1} and \ref{vector2}, respectively. The curves of $\hat{x}_i^1, i=1,\cdots,30$ soon  converge to 0 with some fluctuation. The curves of $\hat{x}_i^2, i=1,\cdots,30$ first decrease and then grow up before the time 40. After a period of time, they converge to a constant, which verify the $\rho$-stability of agents.


\section{Concluding Remarks}

In this paper, we have considered uniform stabilization and social optimality for mean field LQ multiagent systems. 
For finite- and infinite- horizon problems, we design the decentralized control laws by decoupling FBSDEs, respectively, which are further shown to be asymptotically optimal. Some equivalent conditions are further given for uniform stabilization of the systems in different cases. Finally, we compare such decentralized control laws with the asymptotic optimal strategies in previous works.

An interesting generalization is to consider mean field LQ control systems with heterogeneous coefficients by the direct approach \cite{HQL15}.
Also, the variational analysis may be applied to construct decentralized control laws for the nonlinear social control model. 

\appendix
\section{Proof of Theorem \ref{thm3}}\label{app a}
\def\theequation{A.\arabic{equation}}
\setcounter{equation}{0}

To prove Theorem \ref{thm3}, we need a lemma.

\begin{lemma}\label{lem1}
	Let A1) hold. Assume that (\ref{eq8a}) and (\ref{eq9a}) admit a solution, respectively. 
Under the control (\ref{eq15a}), we have
	\begin{equation}
	\max_{0\leq t\leq T}\mathbb{E}\|\hat{x}^{(N)}(t)-\bar{x}(t)\|^2=O({1}/{N}).
	\end{equation}	
\end{lemma}
\emph{Proof.} It follows by (\ref{eq20}) that
\begin{equation*}
\begin{aligned}
d\hat{x}^{(N)}(t)=&\big[(\bar{A}(t)+G)\hat{x}^{(N)}(t)-BR^{-1}B^T(K(t)\bar{x}(t)\cr
&+s(t))+f(t)\big]dt+\frac{1}{N}\sum_{i=1}^N\sigma (t)dW_i(t).
\end{aligned}
\end{equation*}
where $\bar{A}(t)=A-BR^{-1}B^TP(t)$. From (\ref{eq12a}), we have
\begin{equation}\label{eq21}
\begin{aligned}
&\hat{x}^{(N)}(t)-\bar{x}(t)=\Phi(t)[\hat{x}^{(N)}(0)-\bar{x}(0)]\\
&+\frac{1}{N}\sum_{i=1}^N\int_0^t\Phi(t-\tau)\sigma(\tau) dW_i(\tau),
\end{aligned}
\end{equation}
where $\Phi$ satisfies
$\dot{ \Phi}=(\bar{A}(t)+G)\Phi(t),\  \Phi(0)=I.$
By A1), one can obtain
\begin{equation}
\begin{aligned}
\mathbb{E}\big\|  \hat{x}^{(N)}(t)-\bar{x}(t)\big\|^2 \leq&\big\|2\Phi(t)\big\|^2\Big\{\mathbb{E}\big\|\hat{x}^{(N)}(0)-\bar{x}_0\big\|^2\\
&+\frac{1}{N}\int_0^t\big\|
\Phi(\tau)\sigma(\tau)\big\|^2 d\tau\Big\}\cr
\leq \ &  \frac{2}{N}\big\|\Phi(t)\big\|^2\Big\{\max_{1\leq i\leq N}\mathbb{E}\|\hat{x}_{i0}\|^2\\
&+\int_0^T\big\|\Phi(\tau)\sigma(\tau)\big\|^2 d\tau\Big\},
\end{aligned}
\end{equation}
which completes the proof.  $\hfill \Box$

\emph{Proof of Theorem \ref{thm3}}. Note that $\inf_{u\in L^2_{{\mathcal  F}}(0, T; \mathbb{R}^{nr})}J^{\rm F}_{\rm soc}(u)\leq J^{\rm F}_{\rm soc}(\hat{u}). $
To obtain $$\frac{1}{N}J^{\rm F}_{\rm soc}(\hat{u})\leq\frac{1}{N}\inf_{u\in L^2_{{\mathcal  F}}(0, T; \mathbb{R}^{nr})}J^{\rm F}_{\rm soc}(u)+O(\frac{1}{\sqrt{N}}),$$
we only need to prove for any
$u\in \mathcal{U}^{\prime}\stackrel{\Delta}{=}\{u\in L^2_{{\mathcal  F}}(0, T; \mathbb{R}^{nr}): J^{\rm F}_{\rm soc}(u)\leq J^{\rm F}_{\rm soc}(\hat{u})\},$
the following holds:
$$\frac{1}{N}J^{\rm F}_{\rm soc}(\hat{u})\leq\frac{1}{N}J^{\rm F}_{\rm soc}(u)+O(\frac{1}{\sqrt{N}}).$$
We now show that for $u\in \mathcal{U}^{\prime}$,  
$\sum_{i=1}^N\mathbb{E}\int_0^Te^{-\rho t}(\|x_i(t)\|^2+\|u_i(t)\|^2)dt<NC_2$.
{By Lemma \ref{lem2.1}, (P1) is uniformly convex which gives there exists $\delta_0>0$ such that
$$\delta_0 \sum_{i=1}^N\mathbb{E}\int_0^T e^{-\rho t}\|u_i(t)\|^2dt-C\leq J_{\rm soc}^{\rm F}(u).$$
Since $J^{\rm F}_{\rm soc}(\hat{u})<NC_{1}$, we have
$J^{\rm F}_{\rm soc}( {u})<NC_{1}$, 
which implies $ \sum_{i=1}^N\mathbb{E}\int_0^Te^{-\rho t}\|u_i(t)\|^2dt<NC.$} This leads to
$$\mathbb{E}\int_0^T\!e^{-\rho t}\|u^{(N)}\|^2dt\leq\frac{1}{N}\sum_{i=1}^N\mathbb{E}\int_0^T\!e^{-\rho t}\|u_i\|^2dt<C,$$
where $u^{(N)}=\frac{1}{N}\sum_{i=1}^Nu_i.$
By (\ref{eq1}),
\begin{equation*}\label{eq38}
\begin{aligned}
dx^{(N)}(t)=\ & \left[(A+G)x^{(N)}(t)+Bu^{(N)}(t)+f(t)\right]dt\\
&+\frac{1}{N}\sum_{i=1}^N\sigma(t) dW_i(t),
\end{aligned}
\end{equation*}
which 
implies 
$\max_{0\leq t\leq T}\mathbb{E}\|{x}^{(N)}(t)\|^2\leq C.$
Note that
$$x_i(t)=e^{At}x_{i0}+\int_0^te^{A(t-\tau)}[Gx^{(N)}(\tau)+Bu_i(\tau)+f(\tau)]d\tau.$$
We have
\begin{equation}\label{eq22a}
\begin{aligned}
&\sum_{i=1}^N\mathbb{E}\int_0^Te^{-\rho t}\|x_i(t)\|^2dt\\
\leq \ & C\Big(\sum_{i=1}^N\mathbb{E}\|x_{i0}\|^2+N\max_{0\leq t\leq T}\mathbb{E}\|{x}^{(N)}(t)\|^2\\
&+\sum_{i=1}^N\mathbb{E}\int_0^Te^{-\rho t}\|u_i(t)\|^2dt\Big)<NC_2.
\end{aligned}
\end{equation}
By (\ref{eq12a}) and (\ref{eq20}), we obtain that
\begin{equation}\label{eq22b}
\mathbb{E}\int_0^Te^{-\rho t}\big(\|\hat{x}_i(t)\|^2+\|\hat{u}_i(t)\|^2+\|\bar{x}(t)\|^2)dt<C.
\end{equation}

Let $\tilde{x}_i=x_i-\hat{x}_i$,   $\tilde{u}_i=u_i-\hat{u}_i$ and $\tilde{x}^{(N)}=\frac{1}{N}\sum_{i=1}^N \tilde{x}_i$. Then by (\ref{eq1}) and (\ref{eq20}),
\begin{equation}\label{eq32}
d\tilde{x}_i(t)=(A\tilde{x}_i(t)+{G}\tilde{x}^{(N)}(t)+B\tilde{u}_i(t))dt, \  \tilde{x}_i(0)=0.
\end{equation}
From (\ref{eq3}), 
$
J^{\rm F}_{\rm soc}(u)
    =\sum_{i=1}^N(J_i^{\rm F}(\hat{u})+\tilde{J}_i^{\rm F}(\tilde{u})+\mathcal{I}_i)
$,
where
\begin{align*}
\tilde{J}_i^{\rm F}(\tilde{u})\stackrel{\Delta}{=}\mathbb{E}\int_0^T&e^{-\rho t}\big[\|\tilde{x}_i(t)-\Gamma \tilde{x}^{(N)}(t)\|^2_Q+\|\tilde{u}_i(t)\|^2_{R}\big]dt,\\
\mathcal{I}_i=2\mathbb{E}\int_0^T&e^{-\rho t}\Big[\big(\hat{x}_i(t)-\Gamma \hat{x}^{(N)}(t)-\eta(t)\big)^TQ\\
&\times\big(\tilde{x}_i(t)-\Gamma \tilde{x}^{(N)}(t)\big)+\hat{u}_i^T(t)R\tilde{u}_i(t)\Big]dt.
\end{align*}
By Lemma \ref{lem2.1} and
Proposition \ref{prop1}, $\tilde{J}_i^{\rm F}(\tilde{u})\geq 0$. We only need to prove $\frac{1}{N}\sum_{i=1}^N
\mathcal{I}_i=O(\frac{1}{\sqrt{N}})$.
By direct computations, one can obtain 
\begin{equation}\label{eq24a}
\begin{aligned}
\sum_{i=1}^N \mathcal{I}_i
=\ & \sum_{i=1}^N 2\mathbb{E}\int_{0}^{T}e^{-\rho t}\Big\{\tilde{x}_i^T\big[Q(\hat x_i-\Gamma\bar{x}-\eta)\\
&-\Gamma^TQ((I-\Gamma)\bar{x}-\eta)\big]+\sum_{i=1}^N\hat{u}_i^TR\tilde{u}_i\Big\}dt\cr
&+\sum_{i=1}^N2\mathbb{E}\int_{0}^{T}e^{-\rho t}(\hat{x}^{(N)}-\bar{x})^T\Xi \tilde{x}_idt.
\end{aligned}
\end{equation}
By (\ref{eq8a})-(\ref{eq10a}), (\ref{eq32}) and It\^{o}'s formula,
$$  \begin{aligned}
0
=\ &\mathbb{E}\int_{0}^{T}\sum_{i=1}^Ne^{-\rho t}\Big\{- \tilde{x}_i^T\big[ Q\hat{x}_i-Q(\Gamma \bar{x}+\eta)\\ 
&-\Gamma^TQ\left((I-\Gamma) \bar{x}-\eta\right)\big]
-\hat{u}_i^TR\tilde{u}_i)\Big\}dt\\&+N\mathbb{E}\int_{0}^{T}e^{-\rho t}(\hat{x}^{(N)}-\bar{x})^T(G^TP+PG)\tilde{x}^{(N)}dt.
\end{aligned}$$
From this and (\ref{eq24a}), we obtain
$$
\begin{aligned}
\frac{1}{N}\sum_{i=1}^N \mathcal{I}_i=&2\mathbb{E}\int_{0}^{T}e^{-\rho t}(\hat{x}^{(N)}(t)-\bar{x}(t))^T\\
&\times (\Xi +G^TP+PG)\tilde{x}^{(N)}(t)dt.
\end{aligned}
$$
By Lemma \ref{lem1}, (\ref{eq22a}) and (\ref{eq22b}), 
we obtain
$$
\begin{aligned}
\Big|\frac{1}{N}\sum_{i=1}^N \mathcal{I}_i\Big|^2\leq &C\mathbb{E}\int_{0}^{T}e^{-\rho t}\|\hat{x}^{(N)}(t)-\bar{x}(t)\|^2dt\\
&\times \mathbb{E}\int_{0}^{T}e^{-\rho t}\|\tilde{x}^{(N)}(t)\|^2dt,
\end{aligned}
$$
which implies $|\frac{1}{N}\sum_{i=1}^N \mathcal{I}_i|=O(1/\sqrt{N})$.

 Moreover, by (\ref{eq8a}), (\ref{eq11}) and direct calculations,
\begin{align*}
&J^{\rm F}_{\rm soc}(\hat{u})\cr
=&  \sum_{i=1}^N\mathbb{E}\int_0^Te^{-\rho t}\Big(\big\|\hat{x}_i-\Gamma \hat{x}^{(N)}-\eta\big\|^2_{Q}+\|\hat{u}_i\|^2_{R}\Big)dt\cr
=&\sum_{i=1}^N\mathbb{E}\int_0^Te^{-\rho t}\Big(\|\hat{x}_i\|^2_Q+\|\hat{x}^{(N)}\|^2_{\Xi }+\|\eta\|_Q^2\cr
&-2\eta^TQ(I-\Gamma)\hat{x}_i+\|\hat{u}_i\|^2_{R}\Big)dt\cr
=&\sum_{i=1}^N\mathbb{E}\int_0^Te^{-\rho t}\Big(\|\hat{x}_i-\hat{x}^{(N)}\|^2_Q+\|\hat{x}^{(N)}\|^2_{Q-\Xi }+\|\eta\|_Q^2\cr
&-2\eta^TQ(I-\Gamma)\hat{x}_i+\|\hat{u}_i-\hat{u}^{(N)}\|^2_{R}+\|\hat{u}^{(N)}\|^2_{R}\Big)dt\cr
=\!&\sum_{i=1}^N\mathbb{E}\big[\|x_{i0}\!-\!{x}^{(N)}(0)\|^2_P\!+\!\|{x}^{(N)}(0)\|^2_{\Pi}\!+\!2s^T(0)x^{(N)}(0)\big]
\cr
&+\!\sum_{i=1}^N\mathbb{E}\int_0^T\!e^{-\rho t}\Big(\big\|\hat{u}_i\!-\!\hat{u}^{(N)}+R^{-1}B^TP(\hat{x}_i\!-\!\hat{x}^{(N)})\big\|^2_R\cr
&+\big\|\hat{u}^{(N)}+R^{-1}B^T\Pi \hat{x}^{(N)}\big\|^2_{R}\Big)dt+  q_T\cr
= &\sum_{i=1}^N\mathbb{E}\big[\|x_{i0}-{x}^{(N)}(0)\|^2_P+\|{x}^{(N)}(0)\|^2_{\Pi}\cr
&+2s^T(0){x}^{(N)}(0)\big]+  Nq_T+ N\epsilon_T,
\end{align*}
where $q_T$ and $\epsilon_T$ are given by (\ref{eq17a})-(\ref{eq17b}).
$\hfill \Box$

\section{Proofs of Lemma \ref{lem2} and Theorem \ref{thm4}}\label{app b}
\def\theequation{B.\arabic{equation}}
\setcounter{equation}{0}

\emph{Proof of Lemma \ref{lem2}.}
From (\ref{eq21}), we have
$$\begin{aligned}
&\mathbb{E}\int_0^{\infty} e^{-\rho t} \|\hat{x}^{(N)}(t)-\bar{x}(t)\|^2dt \\
\leq\ &  2\mathbb{E}\int_0^{\infty} \left\| e^{(\bar{A}+G-\frac{\rho}{2}I)t}\right\|^2\big\|\hat{x}^{(N)}(0)-\bar{x}(0)\big\|^2dt\cr
&+2\mathbb{E}\int_0^{\infty}e^{-\rho t} \frac{1}{N}  \left\| \int_0^te^{(\bar{A}+G)(t-\tau)}\sigma dW_i(\tau)\right\|^2dt\cr
\leq\ & \frac{2}{N}\int_0^{\infty} \left\| e^{(\bar{A}+G-\frac{\rho}{2}I)t}\right\|^2\mathbb{E}\big\|\max_{1\leq i\leq N}\hat{x}_i(0)\big\|^2dt\cr
&+ \frac{C}{N} \mathbb{E}\int_0^{\infty}e^{-\rho \tau }\|\sigma\|^2\int_\tau^{\infty}\big\|e^{(\bar{A}+\bar{G}-\frac{\rho}{2}I)(t-\tau)}\big\|^2 dtd\tau\cr
\leq\ & O({1}/{N}).
\end{aligned}$$
$\hfill \Box$

\emph{Proof of Theorem \ref{thm4}.} By A1)-A3), Lemmas \ref{lem2a} and \ref{lem2}, we obtain that $\bar{x}\in C_{\rho/2}([0,\infty),\mathbb{R}^n)$ and
$$\mathbb{E}\int_0^{\infty} e^{-\rho t} \left(\big\|\hat{x}^{(N)}(t)-\bar{x}(t)\big\|^2\right)dt=O(\frac{1}{N}),$$
which further gives that 
$\mathbb{E}\int_0^{\infty} e^{-\rho t} \|\hat{x}^{(N)}(t)\|^2dt<\infty.$
Denote $g\stackrel{\Delta}{=}-BR^{-1}B^T((\Pi-P)\bar{x}+s)+Gx^{(N)}+f$. Then
\begin{equation}\label{eq14a}
\hat{x}_i(t)=e^{\bar{A}t}\hat{x}_{i0}+\int_0^te^{\bar{A}(t-\tau)}g(\tau)d\tau+\int_0^te^{\bar{A}(t-\tau)}\sigma dW_i(\tau).
\end{equation}
Note that $\bar{A}-\frac{\rho}{2}I$ is Hurwitz. By Schwarz's inequality,
$$\begin{aligned}
&\mathbb{E}\int_0^{\infty}e^{-\rho t}\|\hat{x}_i(t)\|^2dt\\
\leq &3 \mathbb{E}\int_0^{\infty} \left\| e^{(\bar{A}-\frac{\rho}{2}I)t}\right\|^2\|\hat{x}_{i0}\|^2dt\\
&+3\mathbb{E}\int_0^{\infty}e^{-\rho t}t\int_0^t\Big\|e^{\bar{A}(t-\tau )}g(\tau )\Big\|^2d\tau dt\cr
&+3 \mathbb{E}\int_0^{\infty}e^{-\rho t}\int_0^t tr[e^{\bar{A}^T(t-\tau )}\sigma^T(\tau )\sigma(\tau )e^{\bar{A}(t-\tau )}] d\tau dt\cr
\leq &C+3\mathbb{E}\int_0^{\infty}e^{-\rho \tau }\|g(\tau )\|^2\int_{\tau} ^{\infty}t\big \|e^{(\bar{A}-\frac{\rho}{2}I)(t-\tau )}\big\|^2dtd\tau \cr
&+3 C\mathbb{E}\int_0^{\infty}e^{-\rho \tau }\|\sigma(\tau )\|^2\int_{\tau}^{\infty}\big \|e^{(\bar{A}-\frac{\rho}{2}I)(t-\tau )}\big\|^2dtd\tau \cr
\leq &C_1
\end{aligned}$$
This with (\ref{eq14}) completes the proof.   \hfill $\Box$

\section{Proofs of Theorems \ref{thm5} and \ref{thm6}}\label{app c}
\def\theequation{C.\arabic{equation}}
\setcounter{equation}{0}

\emph{Proof of Theorem \ref{thm5}.}
i)$\Rightarrow$ ii). 
By (\ref{eq20}),
\begin{equation}\label{eq24}
\begin{aligned}
\frac{d\mathbb{E}[\hat{x}_i]}{dt}=\ &\bar{A}\mathbb{E}[\hat{x}_i(t)]-BR^{-1}B^T((\Pi-P)\bar{x}(t)+s(t))\\
&+G\mathbb{E}[\hat{x}^{(N)}(t)]+f(t), \quad \mathbb{E}[\hat{x}_i(0)]=\bar{x}_0.
\end{aligned}
\end{equation}
It follows from A1) that
$\mathbb{E}[\hat{x}_i(t)]=\mathbb{E}[\hat{x}_j(t)]=\mathbb{E}[\hat{x}^{(N)}(t)], \ j\not =i.$
By comparing (\ref{eq18}) and (\ref{eq24}), we obtain $\mathbb{E}[\hat{x}_i(t)]=\bar{x}(t)$. Note $\|\bar{x}(t)\|^2=\big\|\mathbb{E}\hat{x}_i(t)\big\|^2\leq \mathbb{E}\|\hat{x}_i(t)\|^2$.
It follows from (\ref{eq23}) that
\begin{equation}\label{eq25}
\int_0^{\infty}e^{-\rho t} \|\bar{x}(t)\|^2dt<\infty.
\end{equation}
By (\ref{eq18}), we have
$$
\begin{aligned}
\bar{x}(t)=&e^{(A+G-BR^{-1}B^T\Pi)t}\Big[\bar{x}_0\\
&+\int_0^te^{-(A+G-BR^{-1}B^T\Pi)\tau}h(\tau)d\tau\Big] ,
\end{aligned}
$$
where $h=-BR^{-1}B^Ts+f$. By the arbitrariness of $\bar{x}_0$ with (\ref{eq25}) we obtain that $A+G-BR^{-1}B^T\Pi-\frac{\rho}{2}I$ is Hurwitz. That is, $(A+G-\frac{\rho}{2}I, B)$ is stabilizable.
By \cite{AM90}, (\ref{eq16}) admits a unique solution such that $\Pi>0$. Note that $\mathbb{E}[x^{(N)}(t)]^2\leq \frac{1}{N}\sum_{i=1}^N\mathbb{E}[\hat{x}_i^2(t)]$. Then from (\ref{eq23}) we have 
\begin{equation}\label{eq27}
\mathbb{E} \int_0^{\infty}e^{-\rho t}\big\|\hat{x}^{(N)}(t)\big\|^2dt<\infty.
\end{equation}
This leads to $ \mathbb{E} \int_0^{\infty}e^{-\rho t}\|g(t)\|^2dt<\infty$, where $g{=}-BR^{-1}B^T((\Pi-P)\bar{x}+s)+G\hat{x}^{(N)}+f$.
By (\ref{eq14a}), we obtain 
$$
\begin{aligned}
\mathbb{E}\|\hat{x}_i(t)\|^2=\ & \mathbb{E}\left\|e^{\bar{A}t}\left(x_{i0}+\int_0^te^{-\bar{A}\tau}g(\tau)d\tau\right)\right\|^2 \\
&+\mathbb{E}\int_0^ttr\big[\sigma^T(\tau)e^{(\bar{A}^T+\bar{A})(t-\tau)}
\sigma(\tau)\big]d\tau.
\end{aligned}
$$
By (\ref{eq23}) and the arbitrariness of ${x}_{i0}$  we obtain that $\bar{A}-\frac{\rho}{2}I$ is Hurwitz, i.e., $(A-\frac{\rho}{2}I, B)$ is stabilizable.
By \cite{AM90}, (\ref{eq15}) admits a unique solution such that $P>0$.

From (\ref{eq25}) and (\ref{eq27}),
\begin{equation}\label{eq28}
\mathbb{E} \int_0^{\infty}e^{-\rho t} \big\|\hat{x}^{(N)}(t)-\bar{x}(t)\big\|^2 dt<\infty.
\end{equation}
On the other hand, (\ref{eq21}) gives
\begin{equation*}
\begin{aligned}
&\mathbb{E}\big\|\hat{x}^{(N)}(t)-\bar{x}(t)\big\|^2=\mathbb{E}\big\|e^{(\bar{A}+G)t}[\hat{x}^{(N)}(0)-\bar{x}_0]\big\|^2\\
&+\frac{1}{N}\int_0^ttr\big[\sigma^T(\tau)e^{(\bar{A}^T+G^T+\bar{A}+G)(t-\tau)}\sigma(\tau)\big] d\tau.
\end{aligned}
\end{equation*}
By (\ref{eq28}) and the arbitrariness of ${x}_{i0}, i=1,\cdots,N$, we obtain that $\bar{A}+G-\frac{\rho}{2}I$ is Hurwitz.

(ii)$\Rightarrow$(iii).  Let $V(t)=e^{-\rho t}\bar{y}^T(t)\Pi \bar{y}(t)$,
where $\bar{y}$ satisfies
\begin{equation*}
\dot{\bar{y}}(t)=(A+G)\bar{y}(t)+B\bar{u}(t),\quad  \bar{y}(0)=\bar{y}_0.
\end{equation*}
Denote $V$ by $V^*$ when $\bar{u}=\bar{u}^*=-{R^{-1}}B^T\Pi \bar{y}$. By (\ref{eq16}),
\begin{align*}
\frac{dV^*}{dt}=&\bar{y}^T(t)\big[-\rho\Pi+(A+G-B{R^{-1}}B^T\Pi)^T\Pi\\
&+\Pi(A+G-B{R^{-1}}B^T\Pi) \big]\bar{y}(t)\cr
=&\bar{y}^T(t)\big[-(Q-\Xi )-\Pi B{R^{-1}}B^T\Pi\big]\bar{y}(t)\leq 0.
\end{align*}
Note that $V^*\geq0$. Then $\lim_{t\to\infty}V^*(t)$ exists, which implies
\begin{equation}\label{eq43}
\lim_{t_0\to\infty}[V^*(t_0)-V^*(t_0+T)]=0.
\end{equation}

Rewrite $\Pi(t)$ in (\ref{eq11}) by $\Pi_{T}(t)$. Then we have $\Pi_{T+t_0}(t_0)=\Pi_{T}(0)$.
By (\ref{eq11}),
\begin{align*}
&\int_{t_0}^{T+t_0}e^{-\rho t} (\|\bar{y}(t)\|^2_{Q-\Xi }+\|\bar{u}(t)\|^2_R)dt\cr
=&e^{-\rho t_0}\bar{y}^T({t_0})\Pi_{T+t_0}(t_0)\bar{y}({t_0})\cr
&+ \int_0^Te^{-\rho t} \big\|\bar{u}(t)+{R^{-1}}B^T\Pi_{T+t_0}(t_0) \bar{y}(t)\big\|^2_Rdt\cr
\geq&e^{-\rho t_0}\big\|\bar{y}({t_0})\big\|^2_{\Pi_{T+t_0}(t_0)} =e^{-\rho t_0}\big\|\bar{y}({t_0})\big\|^2_{\Pi_{T}(0)}.
\end{align*}
This with (\ref{eq43}) implies
\begin{align*}
&\lim_{t_0\to\infty}e^{-\rho t_0}\big\|\bar{y}({t_0})\big\|^2_{\Pi_{T}(0)}\cr
\leq&\lim_{t_0\to\infty} \int_{t_0}^{T+t_0}e^{-\rho t} (\|\bar{y}(t)\|_{Q-\Xi }^2+\|\bar{u}^*(t)\|^2_R)dt \cr
=&\lim_{t_0\to\infty}[V^*(t_0)-V^*(t_0+T)]=0.
\end{align*}
By A3), one can get that there exists $T>0$ such that $\Pi_{T}(0)>0$ (See e.g. \cite{ZQ16}, \cite{ZZC08}).
Thus, we have $\lim_{t\to\infty}e^{-\rho t}\big\|\bar{y}({t})\big\|^2=0$, which implies that $(A+G-\frac{\rho}{2} I, B)$ is stabilizable. Similarly, we can show $(A-\frac{\rho}{2} I, B)$ is stabilizable.

(iii)$\Rightarrow$(i).
This part has been proved in Theorem \ref{thm4}.  $\hfill \Box$

\emph{Proof of Theorem \ref{thm6}.}
(iii)$\Rightarrow$(i). From \cite{AM90}, (\ref{eq15}) and (\ref{eq16}) admit unique solutions $P\geq0, \Pi\geq0$ such that $A-BR^{-1}B^TP-\frac{\rho}{2}I$ and $A-BR^{-1}B^T\Pi-\frac{\rho}{2}I$
are Hurwitz, respectively. Thus, there exists a unique $s(0)$ such that $s\in C_{\rho/2}([0,\infty),\mathbb{R}^n)$. It is straightforward that $\bar{x}\in C_{\rho/2}([0,\infty),\mathbb{R}^n)$. By the argument in the proof of Theorem \ref{thm4}, (i) follows.
(i)$\Rightarrow$(ii). The proof of this part is similar to that of (i)$\Rightarrow$(ii) in Theorem \ref{thm5}.

(ii)$\Rightarrow$(iii). Since $\Pi\geq0$, there exists an orthogonal $U$ such that
$U^T\Pi U=\left[
\begin{array}{cc}
0 & 0 \\
0& \Pi_{2}
\end{array}
\right],$
where $\Pi_2>0$.
From (\ref{eq15}),
\begin{align}\label{eq40}
\rho U^T\Pi U=&(U^T\bar{\mathbb{A}}U)^TU^T\Pi U+U^T\Pi UU^T\bar{\mathbb{A}}U\cr
&+U^T\bar{Q}  U,
\end{align}
where $\bar{\mathbb{A}}\stackrel{\Delta}{=}A+G-\Pi BR^{-1}B^T\Pi,\bar{Q} =Q-\Xi +\Pi BR^{-1}B^T\Pi $.
Denote
$$U^T\bar{\mathbb{A}}U=\left[\begin{array}{cc}
\bar{\mathbb{A}}_{11}& \bar{\mathbb{A}}_{12} \\
\bar{\mathbb{A}}_{21}& \bar{\mathbb{A}}_{22}
\end{array}\right], \ U^T\bar{Q}  U=\left[\begin{array}{cc}
\Xi _{11}& \Xi _{12} \\
\Xi _{21}& \Xi _{22}
\end{array}\right].$$
By pre- and post-multiplying by $\xi^T$ and $\xi$ where $\xi=[\xi_1^T,0]^T$, it follows that
$$0=\rho\xi^T U^T\Pi U\xi=\xi^T U^T\bar{Q} U\xi.$$
From the arbitrariness of $\xi_1$, we obtain $\bar{Q}  _{11}=0$.
Since $\bar{Q}  $ is semi-positive definite, then $\bar{Q}  _{12}=\bar{Q}  _{21}=0$, and $\bar{Q}  _{22}\geq0$. By comparing each block matrix of both sides of (\ref{eq40}),
we obtain $\bar{\mathbb{A}}_{21}=0$. It follows from (\ref{eq40}) that
\begin{equation}\label{eq40b}
\rho\Pi_2=\Pi_2\bar{\mathbb{A}}_{22}+\bar{\mathbb{A}}_{22}^T\Pi_2+\bar{Q}  _{22}.
\end{equation}

Let $\zeta=[\zeta_1^T,\zeta_2^T]^T=U^T\bar{y}^*$, where $\bar{y}^*$ satisfies $\dot{\bar{y}}^*=\bar{\mathbb{A}}\bar{y}^*$.
Then we have
\begin{align*}
\dot{\zeta_1}&=\bar{\mathbb{A}}_{11}\zeta_1+\bar{\mathbb{A}}_{12}\zeta_2,\cr
\dot{\zeta_2}&=\bar{\mathbb{A}}_{22}\zeta_2.
\end{align*}
By Lemma 4.1 of \cite{W68}, the detectability of $(A+G, (Q-\Xi )^{1/2})$ implies the detectability of $(\bar{\mathbb{A}}, \bar{Q} ^{1/2})$.
Take $\zeta(0)=\xi=[\xi_1^T,0]^T$. Then $\bar{Q} ^{1/2}\bar{y}=\bar{Q} ^{1/2}U\zeta=0$, which together with the detectability of $(\bar{\mathbb{A}}, \bar{Q} ^{1/2})$ implies $\zeta_1\to 0$ and $\bar{\mathbb{A}}_{11}$ is Hurwitz.
Denote $S(t)=e^{-\rho t}\zeta_2^T\Pi_2\zeta_2$. By (\ref{eq40b}),
$$S(T)-S(0)=-\int_0^T\zeta_2(t)^T\bar{Q} _{22}\zeta_2(t)dt\leq 0,$$
which implies $\lim_{t\to\infty}S(t)$ exists.
By a similar argument with the proof of Theorem \ref{thm5}, we obtain $\lim_{t_0\to\infty}e^{-\rho t_0}\big\|\zeta_2(t_0)\big\|^2_{\Pi_{2,T}(0)}=0$
and $\Pi_{2,T}(0)>0$, which gives $\zeta_2\to 0$ and $\bar{\mathbb{A}}_{22}$ is Hurwitz. This with the fact that $\bar{\mathbb{A}}_{11}$ is Hurwitz gives that
$\zeta$ is stable, which leads to (iii). $\hfill \Box$

\end{document}